\documentclass[reqno]{amsart}

\usepackage{graphicx,helvet}
\usepackage{epstopdf}
\usepackage[percent]{overpic}
\usepackage{xcolor}
\usepackage{multirow}
\usepackage{epsfig}
\usepackage{cite}
\usepackage[english]{babel}
\usepackage[utf8]{inputenc}

\def\bcdot{{*}}
\newenvironment{dem}[1][Proof.]{\begin{trivlist}
\item[\hskip \labelsep {\bfseries #1}]}{\end{trivlist}}

\newcommand{\smax}{s^{\rm max}}
\newcommand{\bi}{\mathbf{i}}
\newcommand{\bj}{\mathbf{j}}
\newcommand{\bhalf}{\mathbf{1/2}}
\newcommand{\bef}{\mathbf{e}_1}
\newcommand{\bes}{\mathbf{e}_2}
\newcommand{\bzero}{\mathbf{0}}
\newcommand{\bone}{\mathbf{1}}
\newcommand{\abs}[1]{\left| #1 \right|}

\newcommand{\R}{\mathbb{R}}
\newcommand{\ddim}{s} 

\newtheorem{proposition}{Proposition}[section]

\newtheorem{remark}{Remark}[section]
\newtheorem{lemma}{Lemma}[section]
\newcommand{\bigO}{\mathcal{O}}
\newcommand{\tL}{\textrm{L}}
\newcommand{\tR}{\textrm{R}}
\newcommand{\Z}{\mathbb{Z}}

\numberwithin{equation}{section}

\title[ ACAT methods]{An order-adaptive compact approximation Taylor method for systems of conservation laws}
\author[Carrillo]{H. Carrillo$^{\mathrm{a}}$}
\author[Macca]{E. Macca$^{\mathrm{b}}$}
\author[Russo]{G. Russo$^{\mathrm{c}}$}
\author[Parés]{C. Parés$^{\mathrm{d}}$}
\author[Zorío]{D. Zorío$^{\mathrm{e}}$}

\usepackage[normalem]{ulem}

\begin{document}

\begin{abstract}
		We present a new family of high-order shock-capturing finite difference numerical methods for systems of conservation laws. These methods, called  Adaptive Compact Approximation Taylor (ACAT) schemes, use centered $(2p + 1)$-point stencils, where $p$ may take values in $\{1, 2, \dots, P\}$ according to a new family of smoothness indicators in the stencils.  The methods are based on a combination of a robust first order scheme and the Compact Approximate Taylor (CAT) methods of order  $2p$-order, $p=1,2,\dots, P$ so that they are first order accurate near discontinuities and have order $2p$ in smooth regions, where $(2p +1)$ is the size of the biggest stencil in which  large gradients are not detected. CAT methods, introduced in \cite{CP2019}, are an extension to nonlinear problems of the Lax-Wendroff methods in which the Cauchy-Kovalesky (CK) procedure is circumvented following the strategy introduced in \cite{ZBM2017} that allows one to compute  time derivatives in a recursive way using high-order centered differentiation formulas combined with  Taylor expansions in time.  The expression of ACAT methods for 1D and 2D systems of balance laws are given and the performance is tested in a number of test cases  for several linear and nonlinear systems of conservation laws, including Euler equations for gas dynamics.

\end{abstract} 

\date{\today}


\thanks{$^{\mathrm{A}}$ University of Málaga, Spain.  E-Mail: 
   {\tt hcarrillo@uma.es}}

 \thanks{$^{\mathrm{B}}$ University of Catania, Italy.  E-Mail: 
   {\tt emanuele.macca@phd.unict.it}}

 \thanks{$^{\mathrm{C}}$ University of Málaga, Spain. E-Mail: 
   {\tt pares@uma.es}}

\thanks{$^{\mathrm{D}}$ University of Catania, Italy.  E-Mail: 
   {\tt russo@dmi.unict.it}}

\thanks{$^{\mathrm{E}}$ CI$^{\mathrm{2}}$MA, University of
Concepci\'{o}n, Chile.  E-Mail:  {\tt dzorio@ci2ma.udec.cl}}

\maketitle

\section{Introduction}
Lax-Wendroff methods  for linear systems of conservation laws are based on Taylor expansions in time in which the time derivatives are transformed into spatial derivatives using the governing equations  \cite{LeVeque2007book,Toro2009book,Gideon1971}. The spatial derivatives are then discretised by means of centered high-order differentiation formulas. This procedure allows to derive numerical methods of order $2p$, where $p$ is an arbitrary integer, using a centered  $(2 p + 1)$-points stencil that guaranties the $L^2$ stability, see \cite{CP2019}.

The main difficulty to extend Lax-Wendroff methods to nonlinear problems comes from the transformation of time derivatives into spatial derivatives through the Cauchy-Kovalesky (CK) procedure: this approach may indeed be impractical from the computational point of view because it often requires extended symbolic calculus, ended up into inefficient codes. In the context of ADER methods introduced by Toro and collaborators (see \cite{Ader2001, TitarevToro2002, Schwartzkopff2002}), this difficulty has been circumvented   by replacing the CK procedure by local space-time problems that are solved with a Galerkin method: see \cite{DumbserToro2008}, \cite{PNPM}. 

We follow here the strategy introduced in \cite{ZBM2017} to avoid the CK procedure in which time derivatives are  computed in a recursive way using high-order centered differentiation formulas combined with  Taylor expansions in time. This strategy
leads to high-order Lax-Wendroff Approximated  methods (LAT) that are oscillatory close to discontinuities: in \cite{ZBM2017} they were combined with WENO reconstructions to compute the first time derivatives. The resulting  methods (LAT) give non-oscillatory and accurate results. 

Compact Approximated Taylor methods (CAT) introduced in \cite{CP2019}  circumvent the CK procedure using the same strategy as LAT methods. These methods are compact in the sense that the length of the stencils is minimal: $(2p +1)$-point stencils are used to get order $2p$ compared to $4p +1$-point stencils in LAT methods.  The technique used to reduce the length of the stencil makes the computational cost of a time step in CAT methods bigger than in LAT methods: the Taylor expansions are computed locally, so that the total number of  expansions needed to update the numerical solution is multiplied by $(2p + 1)$.
On the other hand, unlike LAT methods, CAT  methods reduce to the standard high-order Lax-Wendroff methods when applied to linear problems and, due to this, they have better stability properties than LAT and allows one to increase the length of time steps, what compensates the extra cost of every time iteration:
see \cite{CP2019}. 

CAT methods have been also combined with WENO in \cite{CP2019} and \cite{CPZ2020} to avoid oscillations near discontinuities. Nevertheless this combination is not optimal: while the best CAT methods are those of even order, WENO methods have odd accuracy order. Moreover, the restriction on the time step imposed by WENO methods may spoil the advantages of the better stability property of CAT methods. To avoid this, we present in this work a new family of methods in which the oscillations near discontinuities produced by CAT methods are cured by adapting the order of accuracy -- and thus the width of the stencils -- to the smoothness of the solution. To do this,  a new class of smoothness indicators is introduced. 

This paper is organized as follows: Section 2 is devoted to briefly recall LAT and CAT methods. In Section 3, we introduce the Adaptive Compact Approximate Taylor Method (ACAT) and a new family of high order smoothness indicators. In Section 4, the extension to 2D problems of ACAT methods is introduced. In Section 5, the results of the numerical experiments for some selected tests, involving 1D and 2D  linear and nonlinear systems of conservation laws, are given in order to compare  the performance of the ACAT methods with WENO methods. Finally, in Section 6, we draw some conclusions.

\section{Approximate Taylor Methods} \label{sec2}

We consider the one-dimensional system of conservation laws
\begin{equation} \label{cons_law_0}
u_t + f(u)_x = 0 , \quad \quad u(x,0)=u_0(x),  \quad -\infty <x < \infty.
\end{equation}
The solution $u:\R \times \R \rightarrow\R^{\ddim}$ 
is an $\ddim$-dimensional vector of conserved quantities. Taylor expansion in time can be used to update numerical solutions of problem (\ref{cons_law_0}) by  
\begin{equation}\label{atm_0}
u_i^{n+1} = u_i^n +\sum_{k=1}^{m} \frac{\Delta t ^k}{k! }\,\tilde{u}^{(k)}_i + \mathcal{O}\left(\Delta t^{m+1}\right).
\end{equation}
where $\{ x_i \}$ are the nodes of a uniform mesh of step $\Delta x$; $u_i^n$ is a pointwise approximation of the solution at $x_i$ at the time $n \Delta t$, where $\Delta t$ is the time step; and $\tilde{u}^{(k)}_i$ is an approximation of $\partial_t^k u(x_i, n \Delta t)$.

The strategy followed in \cite{ZBM2017} to avoid the CK procedure is based on the equalities  
\begin{equation}\label{u_t^k}
\partial_t^k u = -\partial_x\partial_t^{k-1}f(u).
\end{equation}
that can be easily derived from the equation, if the solutions are assumed to be smooth enough. Numerical approximations of the derivatives appearing at the right-hand side are computed by combining numerical differentiation formulas in space and time with Taylor expansions in a recursive way. 
 For the sake of simplicity, the methods will be only described for the one-dimensional scalar case: extension to systems is straightforward.

\subsection{Lax-Wendroff Approximate Taylor Methods }\label{LAT} 

In Lax-Wendroff Approximate Taylor (LAT) methods,  the time derivatives $\partial_t^{k}u$ are approximated
by applying a first order numerical differentiation formula in space to some approximations 
\begin{equation}\label{f(k-1)_i}
\tilde{f}^{(k-1)}_{i}  \approx \partial_t^{k-1}f(u)(x_{i}, t_n)
\end{equation}
that will be computed by using recursively Taylor expansions in time. 
Here for any function $g(x,t)$ we shall denote by $\tilde{g}_i^{(k)}$ the approximation of the $k$-th time derivative in $(x_i,t_n)$, while 
$g_i^{(k)}$ denotes the corresponding $k$-th space derivative, i.e.
\[
    \tilde{g}_i^{(k)} \cong   \partial^k_tg(x_i,t_n), 
    \quad g_i^{(k)} \cong \partial^k_xg(x_i,t_n).
\]
Centered $(2p +1)$-point numerical differentiation formulas
\begin{equation} \label{F}
f^{(k)}(x_i)   \simeq   D^k_{p,i}(f, \Delta x) =  \frac{1}{\Delta x^k} \sum_{j=-p}^{p} \delta^k_{p,j} f(x_{i+j})
\end{equation}
will be used to compute derivatives. 

The following  notation 
\begin{equation}
 D^k_{p,i}(f_\bcdot, \Delta x)   =   \frac{1}{\Delta x^k} \sum_{j=-p}^{p} \delta^k_{p,j} f_{i+j},
\end{equation}
will be used to indicate that the formula is applied to some approximations $f_i$ of $f$ and not to its exact point values $f(x_i)$. 
In cases where there are two or more indices, the symbol $\bcdot$ will be used to indicate with respect to which the differentiation is applied. For instance:
\begin{eqnarray*}
\partial^k_x u(x_i , t_n) & \simeq & D^{k}_{p,i} (u_\bcdot^n, \Delta x) =  \frac{1}{\Delta x^k} \sum_{j = -p}^p 
\delta^{k}_{p,j} u^n_{i+j}, \\
\partial^k_t u(x_i, t_n) &\simeq& D^{k}_{p,n} (u_i^\bcdot, \Delta t) =  \frac{1}{\Delta t^k} \sum_{r = -p}^p \delta^{k}_{p,r} u^{n+r}_{i}.
\end{eqnarray*}

Once the approximations \eqref{f(k-1)_i} have been computed, the time derivatives of the solution are approximated by:
$$
\partial_t^k u(x_i, t_n) \cong  \tilde u^{(k)}_i =  - D^1_{p,i} (\tilde{f}^{(k-1)}_\bcdot, \Delta x) = - \frac{1}{\Delta x}\ \sum_{j=-p}^p \delta^1_{p,j}  \tilde{f}^{(k-1)}_{i+j}.
$$
where $p$ is adequately chosen so that the local discretization error is of order $O(\Delta x^{m +1})$.

Following a recursive procedure, the approximation of the time derivatives are used to compute approximations of the flux forward and backward in time using Taylor expansions. Once all the time derivatives are approximated,  \eqref{atm_0} is used to update the numerical solutions.

The procedure can be summarized as follows:

\begin{enumerate}

\item Define
$$
\tilde{f}^{(0)}_i = f(u^n_{i}).
$$

\item Compute
 \begin{equation}\label{lat_ut}
 \tilde{u}^{(1)}_{i} = - D^1_{p,i}(\tilde{f}^{(0)}_\bcdot, \Delta x).
 \end{equation}

\item  For $k = 2, \dots, m$:

\begin{enumerate}

\item Compute  an approximation of $f(u(x_i,t_{n+r}))$ as  
$$
\tilde{f}^{k-1,n+r}_{i} = f \left(  u^n_{i} + \sum_{l=1}^{k-1} \frac{(r \Delta t)^l}{l!} \tilde{u}^{(l)}_{i} \right), \quad 
r = -p, \dots, p.
$$
where the approximate Taylor expansion of the function  $u(x_i,t_{n+r})$ has been used.

\item Compute 
\begin{equation}\label{derf(k-1)}
\tilde{f}^{(k-1)}_{i} =  D^{k-1}_p(\tilde{f}^{k-1, \bcdot}_i, \Delta t).
\end{equation}

\item Compute 
\begin{equation}\label{deru(k)}
    \tilde{u}  ^{(k)}_{i} = - D^1_{p,i}(\tilde{f}^{(k-1)}_\bcdot, \Delta x).
\end{equation}

\end{enumerate}

\item Update the solution by (\ref{atm_0}).

\end{enumerate}

The order of the method is $\min(m, 2p)$.

\begin{remark}
 Although, for the sake of clarity, we present $m$ and $p$ here like two arbitrary independent positive integers, in \cite{ZBM2017} 
 $m$ is an odd number (since the method is combined with WENO reconstructions) and $p$ is chosen adequately to obtain order $m$. More precisely, in formulas \eqref{deru(k)},
 $$
 p = \left\lceil  \frac{m+1-k}{2}  \right\rceil,
 $$
 where $\lceil \cdot \rceil$ is the ceiling function, and in formulas \eqref{derf(k-1)}
 $$
 p = \frac{m- 1}{2}.
 $$
 \end{remark}
This family of methods can be also written in conservative form. To see this, let us introduce the family of interpolatory numerical differentiation formulas
\begin{equation}\label{upwF}
f^{(k)}(x_i + q \Delta x) \simeq A^{k,q}_{p,i}(f, \Delta x) = \frac{1}{\Delta x^k} \sum_{j = -p + 1}^p \gamma^{k,q}_{p,j} f(x_{i+j}),
\end{equation}
 that approximates the $k$-th derivative of a function at the point $x_i + q \Delta x$ using its values at the $2p$  points $x_{i-p+1}, \dots, x_{i+p}$. The symbol $\bcdot$ will be used again to indicate with respect to which index the differentiation is performed. 

\begin{remark}
The coefficients $\delta^k_{p,j}$ and $\gamma^{k,q}_{p,j}$  of the differentiation formulas can be recursively computed using the algorithm introduced in \cite{Fornberg1}. See also \cite{CP2019}. 
\end{remark}

The following relation holds (see \cite{CP2019}):
\begin{equation}\label{F2}
D^k_{p,i}(f, \Delta x)   =  \frac{1}{\Delta x} \left(A^{k-1, 1/2}_{p,i}(f, \Delta x) -  A^{k-1, 1/2}_{p,i-1}(f, \Delta x)\right).
\end{equation}
Using this equality with $k = 1$, the methods can be written in the form
  \begin{equation}\label{cons}
u_i^{n+1} = u_i^n + \frac{\Delta t}{\Delta x}\left( F^p_{i-1/2} - F^p_{i+1/2}\right),
\end{equation}
where
\begin{equation}\label{cat1}
F^p_{i+1/2}  = \sum_{k=1}^{m} \frac{\Delta t^{k-1}}{k!}A^{0, 1/2}_{p,i}(\tilde{f}_{\bcdot}^{(k-1)}, \Delta x).  
\end{equation}

\subsection{Compact Approximate Taylor methods}

CAT methods were designed in \cite{CP2019} as a variant of the  previous methods that properly generalize  the Lax-Wendroff methods for linear systems. These methods are based on the conservative expression (\ref{cons},\ref{cat1}) but the difference is that now the numerical flux $F_{i+1/2}$ is computed using only the values 
\begin{equation*} 
u^n_{i -p +1}, \dots, u^n_{i+ p}, 
\end{equation*}
so that $u_i^{n+1}$ is updated using only the values at the centered $(2p + 1)$-point stencil. The numerical flux is given by
\begin{equation}\label{cat2}
F^p_{i+1/2}  = \sum_{k=1}^{m} \frac{\Delta t^{k-1}}{k!}A^{0, 1/2}_{p,0}(\tilde{f}_{i,\bcdot}^{(k-1)}, \Delta x),  
\end{equation}
where 
\begin{equation}\label{f(k-1)_ij}
\tilde{f}^{(k-1)}_{i,j}  \approx \partial_t^{k-1}f(u)(x_{i+j}, t_n), \quad j=-p+1,\dots, p
\end{equation}
are \textit{local} approximations of the time derivatives of the flux.

\begin{remark}
The formula $A^{0, 1/2}_{p,0}$ appearing in \eqref{cat2} indicates the Lagrange interpolation (since the order of differentiation is zero) of 
$$
\tilde{f}_{i,j}^{(k-1)}, \quad j= -p+1, \dots, p,
$$
evaluated at $x_{i+1/2}$, whose 'local index' is $j = 0 + 1/2$: that is the reason of the subindex 0 and  the superindex $1/2$.
\end{remark}
By \textit{local} we mean that these approximations depend on the stencil, i.e.
$$
i_1 + j_1 = i_2 + j_2  \not \Rightarrow \tilde{f}^{(k-1)}_{i_1,j_1} = \tilde{f}^{(k-1)}_{i_2,j_2}.
$$
Local approximations of the time derivatives of the solution
\begin{equation*}
\tilde{u}  ^{(k)}_{i,j} \cong \partial_t^{(k)}u(x_{i+j}, t_n), \quad j=-p+1,\dots, p
\end{equation*}
are then computed by using the non-centered differentiation formulas applied to the discrete version of Eq. \eqref{f(k-1)_i}
\begin{equation*}
\tilde{u} ^{(k)}_{i,j} = - A^{1,j}_{p,0}(\tilde{f}^{(k-1)}_{i, \bcdot}, \Delta x)
= - \frac{1}{\Delta x} \sum_{r=-p +1}^p \gamma^{1,j}_{p,r} \tilde{f}^{(k-1)}_{i, r}.
\end{equation*}
These approximations of the time derivatives are then used to compute the approximations of the flux forward and backward in time using Taylor expansions in a recursive way. 
The procedure to compute $F^p_{i+1/2}$ for the node $i$ is given as follows:

\begin{enumerate}

\item  {Define}
$$
\tilde{f}^{(0)}_{i,j}=f(u^n_{i+j}), \quad  j = -p+1, \dots, p.
$$

\item  {For $k = 2 \dots m$:}

\begin{enumerate}

\item Compute
\begin{equation*}
 \tilde{u}^{(k-1)}_{i,j} = - A^{1,j}_{p,0}(\tilde{f}^{(k-2)}_{i,\bcdot}, \Delta x). 
\end{equation*}

\item Compute 
$$
\tilde{f}^{k-1,n+r}_{i,j} = f \left(  u^n_{i+j} + \sum_{l=1}^{k-1} \frac{(r \Delta t)^l}{l!} \tilde{u}^{(l)}_{i,j} \right), \quad  j, r = -p+1, \dots, p.
$$

\item Compute
$$
\tilde{f}^{(k-1)}_{i,j} =   A^{k-1,0}_{p,n}( \tilde{f}^{k-1, \bcdot}_{i,j}, \Delta t),\quad  j = -p+1, \dots, p.
$$

\end{enumerate}

\item Compute $F^p_{i+1/2}$ by \eqref{cat2}

\end{enumerate}

Once the numerical fluxes have been computed, the numerical solution is updated by using \eqref{cons}.

In \cite{CP2019} it has been shown that the order of the method is $\min(m, 2p)$ so that the optimal choice is $m = 2p$; the corresponding numerical method named CAT$2p$ reduces  to the $2p$-order version of standard Lax-Wendroff method for linear problems.
CAT$2p$ is linearly stable under the standard CFL-1 condition (see \cite{CP2019}).

It can be easily checked that the numerical flux of CAT2 writes as follows:
\begin{align}\label{CAT2_flux} 
F^1_{i+1/2} = \frac{1}{4}(\tilde{f}^{1,n+1}_{i,1}+ \tilde{f}^{1,n+1}_{i,0} + f^{n}_{i+1} + f^{n}_{i}), 
\end{align}
where
\begin{equation} \label{CAT2_flux_2}
\tilde{f}^{1,n+1}_{i,j} = f\left( u^n_{i+j} - \frac{\Delta t}{\Delta x}\bigl( f(u^n_{i+1})-f(u^n_{i})\bigr)\right), \quad j=\{0,1\}.
\end{equation}
This numerical flux reduces to the standard 
Lax-Wendroff second order numerical flux for $f(u) = au$. 
The explicit form of the numerical flux of CAT4 can be found in \cite{CP2019}.

\section{Adaptive Compact Approximate Taylor Method} 

 Although Compact Approximate Taylor methods are linearly  stable in the $L^2$ sense under the usual CFL-1 condition, they may  produce strong oscillations close to a discontinuity of the solution.  Two different techniques were considered in  \cite{CP2019} to avoid these oscillations: to combine CAT2 with a first order robust method using a flux limiter (FL-CAT$2$ method) or, following \cite{ZBM2017}, to use WENO reconstructions to compute the first order time derivatives (WENO-CAT methods). See also \cite{CPZ2020}. 
 
 Here we follow a different strategy and select automatically the stencil used to compute $F_{i+1/2}$ so that its length is maximal among those for which the solution is smooth. More specifically, let us suppose that solutions $\{ u_i^n \}$ at time $n\Delta t$  have been computed. The maximum length of the stencil to compute $F_{i+1/2}$ is set to, say, $2 P$, where $P$ is a natural number. Then, the candidate stencils to compute $F_{i+1/2}$ are
 $$
 S_p = \{x_{i-p +1}, \dots, x_{i+p}\}, \quad p = 1, \dots, P.
 $$
 In order to select the stencil, some smoothness indicators $\psi^p_{i+1/2}$, $p = 1, \dots, P$ are computed such that:
\begin{equation}\label{condiconesS}
\psi^p_{i+1/2} \approx \left\{
\begin{array}{cl}
1 & \mbox { if $\{u_i^n\}$ is 'smooth' in $S_p $,}\\
0 & \mbox {otherwise.} 
\end{array}\right.
\end{equation}
Define now:
$$
\mathcal{A} = \{  p \in \{1, \dots, P \} \ s.t.\ \psi^p_{i+1/2} \cong 1\}.
$$
The idea would be then to define:
$$
F^A_{i+1/2} = \begin{cases}
F^{lo}_{i+1/2} & \text{if $\mathcal{A} = \emptyset$;}\\
F^{p_s}_{i+1/2} & \text{where $p_s = \max(\mathcal{A})$ otherwise;} \\
\end{cases}
$$

where $F^{p_s}_{i+1/2}$ is the numerical flux of CAT$2p_s$ and $F^{lo}_{i+1/2}$ is a robust first order numerical flux.
Nevertheless,  it is not possible to determine if the solution is smooth or not in the stencil
$S_1$ where only two values $u^n_i$, $u^n_{i+1}$ are available. Therefore, what will be done in practice is to define:
\begin{equation}\label{admissibleset}
\mathcal{A} = \{  p \in \{2, \dots, P \} \ s.t.\ \psi^p_{i+1/2} \cong 1\}.
\end{equation}
and then:
\begin{equation}\label{numfluxacat}
F^A_{i+1/2} = \begin{cases}
F^{*}_{i+1/2} & \text{if $\mathcal{A} = \emptyset$;}\\
F^{p_s}_{i+1/2} & \text{where $p_s = \max(\mathcal{A})$ otherwise;} \\
\end{cases}
\end{equation}
where $F^*_{i+1/2}$ is the numerical flux of the FL-CAT$2$ (that uses the stencil $S_2$ as well). 
In what follows, we recall first the expression of the FL-CAT$2$  numerical flux;  next, we introduce the smoothness indicators; then, we summarize the expression of the high-order ACAT methods; and finally we briefly discuss its application to systems of conservation laws.

\subsection{FL-CAT2 numerical flux}
Let us consider the scalar conservation law \eqref{cons_law_0} with  $m = 1$.
The expression of the FL-CAT2 numerical flux is as follows:
\begin{align}\label{ACAT2_flux} 
F^*_{i+1/2} = \psi^1_{i+1/2} \, F^1_{i+1/2} +(1-\psi^1_{i+1/2}) \, F^{lo}_{i+1/2},  
\end{align}
where $F^1_{i+1/2}$ is given by \eqref{CAT2_flux}-\eqref{CAT2_flux_2};  $F^{lo}_{i+1/2}$ is a first-order robust numerical flux; and $\psi^1_{i+1/2}$ is a standard flux limiter: 
\begin{equation}\label{indicador1}
\psi_{i+1/2}^1  =  \psi^1(r_{i+1/2}), 
\end{equation}
where
\begin{equation}\label{wavefuction}
r_{i+1/2}=\frac{\Delta upw}{\Delta loc}
= \left\{
\begin{array}{cl}
\displaystyle r^-_{i+1/2}:= \frac{u^n_i-u^n_{i-1}}{u^n_{i+1}-u^n_{i}} & \mbox {if } a_{i+1/2}  >0, \\
\displaystyle  r^+_{i+1/2}:= \frac{u^n_{i+2}-u^n_{i+1}}{u^n_{i+1}-u^n_{i}} & \mbox {if } a_{i+1/2}  \leq 0;
\end{array}\right.
\end{equation}
and  $a_{i+1/2}$ is an estimate of the wave speed like for instance Roe's intermediate speed:
$$a_{i+1/2} = \begin{cases}  \displaystyle \frac{f(u^n_{i+1}) - f(u^n_i)}{u^n_{i+1} - u^n_i} & \text{ if $u^n_i \not= u^n_{i+1};$}\\
f'(u^n_i) & \text{otherwise.}
\end{cases}$$  
An alternative that avoids the computation of an intermediate speed was introduced in \cite{Toro2009book}: it consists in defining
\begin{equation}\label{local_indicators}
\psi^1_{i+1/2}=\min( \psi^1(r_{i+1/2}^+), \psi^1(r^-_{i+1/2})).
\end{equation}
This expression of the smoothness indicator is especially useful for systems: see Section \ref{ss:systems}.

\subsection{Smoothness indicators}	\label{ss:smoothness}

Let us introduce a new family of local smoothness indicators $\psi^{p}_{i+1/2}$, $p\geq2$, for scalar conservation laws and analyze their properties. 

Given the nodal approximations $f_i$ of a function $f$ at the stencil $S_p$, $p\geq2$, centered at $x_{i+1/2}$, first define the lateral weights:
\begin{equation}
I_{p,L}:=\sum_{j=-p+1}^{-1}(f_{i+1+j}-f_{i+j})^2+\varepsilon,\quad
I_{p,R}:=\sum_{j=1}^{p-1}(f_{i+1+j}-f_{i+j})^2+\varepsilon,
\end{equation}	
where $\varepsilon$ is a small quantity that is added to prevent the lateral weights to vanish when the function is constant. Next, compute:
\begin{equation}
I_{p}:=\frac{I_{p,L}I_{p,R}}{I_{p,L}+I_{p,R}}.
\end{equation}
Finally, define the smoothness indicator of the stencil of $S_p$ by 
\begin{equation}\label{indicadores_locales}
\psi_{i+1/2}^p:=\left( \frac{I_p}{I_p+\tau_p}\right),
\end{equation}
where 
\begin{align}
\tau_p:=& \left( \Delta^{2p-1}_{i-p +1} f\right)^2.
\end{align}
Here, $ \Delta^{2p-1}_{i-p +1} f$ represents the undivided difference of 
$\{ f_{i-p+1}, \dots, f_{i+p}\}$:
\begin{align}
\Delta^{2p-1}_{i-p +1} f =&  (2p - 1)! \sum^{p}_{j=-p+1}\,\gamma^{2p-1, 1/2}_{p,j} \, f^n_{i+j}.
\end{align}
Before going into technical details, let us give a motivation of this choice. If data in the stencil $S_p$ are smooth, then
$$
I_{p,L} = O(\Delta x^2), \quad I_{p,R} = O(\Delta x^2), \quad \tau_p = O(\Delta x^{4p}).
$$
Since
$$
\frac{1}{I_p} = \frac{1}{I_{p,L}} + \frac{1}{I_{p,R}}
$$
then $I_p = O(\Delta x^2)$ and thus
$$
\psi_{i+1/2}^p =  \frac{I_p}{I_p+\tau_p} =  \frac{O(\Delta x^2)}{O(\Delta x^2)+ O(\Delta x^{4p})} \approx 1.
$$
On the other hand, if there is an isolated discontinuity in the stencil then 
$$
\tau_p = O(1)
$$
and 
$$
I_{p,L} = O(1), \quad I_{p,R} = O(\Delta x^2)
$$
or 
$$
I_{p,L} = O(\Delta x), \quad I_{p,R} = O(1).
$$
In both cases $I_p = O(\Delta x^2)$ and thus:
$$
\psi_{i+1/2}^p =  \frac{I_p}{I_p+\tau_p} =  \frac{O(\Delta x^2)}{O(\Delta x^2)+ O(1)} \approx 0.
$$
Nevertheless, in the case of smooth data, special care has to be taken if there is a critical point in the stencil, since in this case the order of $I_p$ depends on the order of the critical point, what can prevent the smoothness indicator to be close of 1, as it will be seen in Propositions \ref{prop:accuracy}-\ref{prop:accuracyp=2mod} below. The following definition is assumed in these results: a point $x$ is said to be a critical point of $f$ of order $n$ if $f^{(j)}(x) = 0$, $j = 1, \dots, n$ and $f^{(n+1)}\not= 0$.

Before analysing the smoothness indicators, let us introduce some definitions and notation, taken from 
\cite{ZBBM2019}: we refer to Section 2.1 of this reference for further details. 

Given $\alpha\in\mathbb{R}^+$ and $f:(0, h^*) \mapsto \mathbb{R}$ with $h^* \in (0, \infty]$, the notation  
$f(h)=\mathcal{O}(h^\alpha)$ means, as usual, that
$$ \limsup_{h\rightarrow0^+}\left| \frac{f(h)}{h^\alpha} \right| <+\infty,$$ 
and the notation $f(h)=\mathcal{\bar{O}}(h^\alpha)$ means that
$$ \limsup_{h\rightarrow0^+}\left| \frac{f(h)}{h^\alpha} \right| <+\infty \quad \mathrm{and} \quad \liminf_{h\rightarrow0^+}\left| \frac{f(h)}{h^\alpha} \right|>0. $$

If $f, g: (0, h^*) \mapsto \mathbb{R}$ and $\alpha, \beta$ are two positive real numbers, the following relations hold:
\begin{eqnarray*}
& & f(h) =  \mathcal{O}(h^\alpha), \quad g(h) = \mathcal{O}(h^\beta) \implies f(h)g(h)  = \mathcal{O}(h^{\alpha + \beta}); \\
& & f(h) = \mathcal{\bar{O}}(h^\alpha), \quad g(h) = \mathcal{\bar{O}}(h^\beta) \implies f(h)g(h) =  \mathcal{\bar{O}}(h^{\alpha+\beta});\\
& & f > 0, f(h) = \mathcal{\bar{O}}(h^\alpha) \implies f(h)^{-1} =  \mathcal{\bar{O}}(h^{1/\alpha}).
\end{eqnarray*}

\begin{lemma} Let $c, d, z \in \mathbb{R}$. Assume that
$$
\left\{ \begin{array}{ll}
f^{(j)}(z) = 0 \text{ for } j= 1, \dots, k, \quad f^{(k+1)}(z) \not= 0, \text{ and } f\in \mathcal{C}^{k+2} & \text{if $c + d\not= 0$;}\\
f^{(2j-1)}(z) = 0 \text{ for } j= 1, \dots, n, \quad f^{(2n+1)}(z) \not= 0, \text{ and } f\in \mathcal{C}^{2n+2} & \text{if $c + d = 0$.}
\end{array}
\right.
$$
Then 
$$f(z + dh) - f(z-dh) = \mathcal{\bar{O}}(h^s),$$
where
$$
s = \left\{ \begin{array}{ll}
k + 1 & \text{if $c + d\not= 0$;}\\
2n + 1 & \text{if $c + d = 0$.}
\end{array}
\right.
$$

\end{lemma}

From this lemma, whose proof is given in \cite{ZBBM2019}, one can deduce that, given the values
$f_j = f(x_j)$, $j = i-p+1, \dots, i+p$ of a smooth enough function  $f$ in the stencil $S_p$, 
the following estimates hold:
  $$ f_{j+1} - f_j = \mathcal{O}(h), \quad j = i-p+1, \dots, i + p -1$$
  if the stencil does not contain any critical point of $f$;
 \begin{equation}\label{localorder}
  f_{j+1} - f_j = \mathcal{\bar{O}}(h^{k+1}),  \quad j = i-p+1, \dots, i + p -1,
 \end{equation}
if the stencil contains a critical point $x^*$ of even order $k$ or a critical point of odd order
that is not located at the center of any sub-interval of the stencil. 

Finally, if there exists
$i_0$ such that $x^* = 0.5(x_{i_0} + x_{i_0 +1})$ is a critical point of odd order, then
\eqref{localorder} holds for every $j \not= i_0$ and
 \begin{equation}\label{malasuerte}
  f_{i_0+1} - f_{i_0} = \mathcal{\bar{O}}(h^{2n + 1}) 
 \end{equation}
where $2n + 1$ is the first odd number such that 
$$f^{(2n+1)}(x^*) \not= 0.$$

Let us analyze the behavior of the smoothness indicators \eqref{indicadores_locales} assuming that $\varepsilon = 0$ (the role of $\varepsilon$ is only relevant for the implementation of the method):

\begin{proposition}\label{prop:accuracy}
Let $f_j = f(x_j)$, $j = i-p+1, \dots, i+p$ be the values of a function  $f$ in the stencil $S_p$, with $p>2$.
The following estimates hold:
	\[{
	\psi^p_{i+1/2}=\begin{cases}
		1 - \bigO(\Delta x^{4(p-1)-2k}) & \textnormal{if $f\in \mathcal{C}^{\max(2p-1,k+2)}$;}\\
	\bar\bigO(\Delta x^{2(k+1)}) & \textnormal{if  $f$ is piecewise $\mathcal{C}^{k+2}$ and $S_p$ contains an isolated jump discontinuity of $f$;}	\end{cases}}
	\]
	where $k = 0$ if there is no critical point of $f$ in $S_p$ or $k$ equal to the order of the critical point if there is one. 

\end{proposition}
\begin{dem}
If $f \in C^{2p-1}$ there exists $\xi$ such that
$$
\Delta^{2p-1}_{i-p +1} f = (2p -1)! f^{(2p -1)}(\xi) \Delta x^{2p -1},
$$
and thus
$$
\Delta^{2p-1}_{i-p +1} f  = \bigO(\Delta x^{2p-1}),
$$
what implies
	$$\tau_p=\bigO(\Delta x^{4p-2}).$$
On the other hand, if $S_p$ contains an isolated jump discontinuity, then
$$
\Delta^{2p-1}_{i-p +1} f  = \bigO(1),
$$
and thus	
	$$\tau_p=\bar\bigO(1).$$
	From the discussion above, the estimate
	$$
	f_{j+1}- f_j = \bar\bigO(\Delta x^{k + 1}),
	$$
holds for every $j \in i-p+1, \dots, i+p-1$ with the exception of at most one index $i_0$, in which the order is higher. 

Nevertheless, since  both $I_{p,L}$ and $I_{p,R}$ are the sum of at least two terms of the form $(f_{j+1}-f_j)^2$, we can conclude that 
	$$I_{p,L}=\bar\bigO(\Delta x^{2+2k}),\quad
	I_{p,R}=\bar\bigO(\Delta x^{2+2k}).$$
	Hence:
	$$I_p =\frac{I_{p,L}I_{p,R}}{I_{p,L}+I_{p,R}}=\frac{\bar\bigO(\Delta x^{2+2k})\bar\bigO(\Delta x^{2+2k})}{\bar\bigO(\Delta x^{2+2k})+\bar\bigO(\Delta x^{2+2k})}=\frac{\bar\bigO(\Delta x^{4+4k})}{\bar\bigO(\Delta x^{2+2k})}=\bar\bigO(\Delta x^{2+2k}).$$
	Now, if $S_p$ contains a discontinuity, then, by construction, there exists a side $\alpha \in\{\tL,\tR\}$ 
	such that 	$I_{p,\alpha}=\bar\bigO(1)$ (the side that contains the discontinuity) while the other side, $\beta \in\{\tL,\tR\}\setminus\{\alpha\}$, satisfies $I_{p,\beta}=\bar\bigO(\Delta x^{2+2k})$.	Therefore
	$$I_p=\frac{I_{p,L}I_{p,R}}{I_{p,L}+I_{p,R}}=\frac{I_{p,\alpha}I_{p,\beta}}{I_{p,\alpha}+I_{p,\beta}}=\frac{\bar\bigO(1)\bar\bigO(\Delta x^{2+2k})}{\bar\bigO(1)+\bar\bigO(\Delta x^{2+2k})}=\frac{\bar\bigO(\Delta x^{2+2k})}{\bar\bigO(1)}=\bar\bigO(\Delta x^{2+2k}).$$
	Combining the above results, we have that, if $f$ is smooth:
	\begin{align*}
	\psi^p_{i+1/2}=&\frac{I_p}{I_p+\tau_p}=\frac{1}{\displaystyle1+\frac{\tau_p}{I_p}}=\frac{1}{\displaystyle1+\frac{\bigO(\Delta x^{4p-2})}{\bar\bigO(\Delta x^{2+2k})}}=\frac{1}{\displaystyle1+\bigO(\Delta x^{4(p-1)-2k})} = 1-\bigO(\Delta x^{4(p-1)-2k}).
	\end{align*}
	On the other hand, if $S_p$ contains a discontinuity, then
	\begin{align*}
	\psi^p_{i+1/2}=&\frac{I_p}{I_p+\tau_p}=\frac{1}{\displaystyle1+\frac{\tau_p}{I_p}}=\frac{1}{\displaystyle1+\frac{\bar\bigO(1)}{\bar\bigO(\Delta x^{2+2k})}}=\frac{1}{1+\bar\bigO(\Delta x^{-2(k+1)})}=\bar\bigO(\Delta x^{2(k+1)}),
	\end{align*}
	which finishes the proof.
	\hfill$\Box$
\end{dem}

Observe that the indicator 	$\psi^p_{i+1/2}$ is able to detect smoothness in the presence of a critical point whose
order is lower than  $2(p-1)$.

In the case $p = 2$ similar arguments lead to prove the following estimates:

\begin{proposition}\label{prop:accuracyp=2}
Let $f_j = f(x_j)$, $j = i-1, \dots, i+2$ be the values of a function  $f$ in the stencil $S_2$.
The following estimates hold:
	\[{
	\psi^2_{i+1/2}=\begin{cases}
		1 - \bigO(\Delta x^{4-2k}) & \textnormal{if $f\in \mathcal{C}^3$;}\\
	\bar\bigO(\Delta x^{2(k+1)}) & \textnormal{if  $f$ is piecewise $\mathcal{C}^{k+2}$ and $S_p$ contains an isolated jump discontinuity of $f$;}	\end{cases} }
	\]
	where $k = 0$ if there is no critical point of $f$ in $S_2$ and $k =1$ if there is a critical point $x^*$ of order 1 such that 
	$f^{(3)}(x^*) \not= 0$ or such that $x^* \not= 0.5(x_{j} + x_{j+1})$ for $j= i-1,  i+1$. 
	
	\end{proposition}
Nevertheless, the estimate cannot be proved when $S_2$ includes a critical point of order 1 located at $0.5(x_{i-1} + x_i)$ or 
$0.5(x_{i+1} + x_{i+2})$ and such
that $f^{(3)}(x^*) \not= 0$: the argument in the proof of Proposition \ref{prop:accuracy} cannot be used since there is only one term in the definition of the local weights.   This is not a limitation in many applications, since this situation is very specific and, even if it happens, unless there is a discontinuity close to the critical  point, smoothness will be detected by at least one of the indicators 	$\psi^p_{i+1/2}$ with
$p > 2$ so that the stencil $S_p$ will be used to update the solution. In any case, the smoothness indicator for $p = 2$ can be modified to 
properly handle these situations as follows: compute the couple of lateral weights:
\begin{eqnarray}
I^1_{2,L} & := & (f_{i}-f_{i-1})^2+\varepsilon,\quad I^1_{2,R}:=  (f_{i+1} - f_i)^2 +  (f_{i+2}- f_{i+1})^2+\varepsilon,\\
I^2_{2,L} & := & (f_{i}-f_{i-1})^2+(f_{i+1} - f_i)^2 + \varepsilon,\quad I^2_{2,R}:= (f_{i+2}- f_{i+1})^2+\varepsilon.
\end{eqnarray}	
Next, compute:
\begin{equation}
I^j_{2}:=\frac{I^j_{2,L}I^j_{2,R}}{I^j_{2,L}+I^j_{2,R}}, \quad j= 1, 2.
\end{equation}
and then, the smoothness indicator of the stencil $S_2$  is given by 
\begin{equation}\label{indicadorlocal2mod}
\widetilde{\psi}_{i+1/2}^2:=\max\left( \frac{I^1_2}{I^1_2+\tau_2}, \frac{I^2_2}{I^2_2+\tau_2}\right).
\end{equation}
The following estimate can be then proved:
\begin{proposition}\label{prop:accuracyp=2mod}
Let $f_j = f(x_j)$, $j = i-1, \dots, i+2$ be the values of a function  $f$ in the stencil $S_2$.
The following estimates hold:
	\[
	{
	\widetilde{\psi}^2_{i+1/2}=\begin{cases}
		1 - \bigO(\Delta x^{4-2k}) & \textnormal{if $f\in \mathcal{C}^3$;}\\
	\bar\bigO(\Delta x^{2(k+1)}) & \textnormal{if  $f$ is piecewise $\mathcal{C}^{k+2}$ and $S_p$ contains an isolated jump discontinuity of $f$;}	\end{cases}}
	\]
	where $k = 0$ if there is no critical points of $f$ in $S_2$ or $k =1$ if there is a critical point $x^*$ or order 1. 
	\end{proposition}
\begin{dem}
    The arguments of the proof of Proposition \eqref{prop:accuracy} are used again. The difference comes from the case in which there is
    a critical point of order 1 located at  at $0.5(x_{i-1} + x_i)$ or 
$0.5(x_{i+1} + x_{i+2})$ and such
that $f^{(3)}(x^*) = 0$. In this case, there exists $j \in \{ 1, 2 \}$ (the one in which the sub-interval with the critical point and the
central sub-interval are considered together in the same lateral weight) such that
$$
 \frac{I^j_2}{I^j _2+\tau_2} = 1 - \bigO(\Delta x^{2}).
 $$
Using this estimate the proof is concluded as in Proposition \eqref{prop:accuracy} 
\end{dem}
 
Let us remark finally that the smoothness indicators \eqref{indicadores_locales} and \eqref{indicadorlocal2mod} have finally the following homothetic invariance property: given a function $f$ and positive numbers $\alpha$, $\beta$, define
$$
g(x) = \alpha f(\beta x).
$$
Then the smoothness indicator of $f$ at a stencil $S_p$ centered at $x_{i+1/2}$ in a mesh with step $\Delta x$ is equal to the smoothness indicator of $g$ at the stencil $S_p$ centered at $\beta x_{i+1/2}$ in a mesh with step $\beta \Delta x$. This property is very important in practice to have smoothness indicators whose behaviour do not depend on $\Delta x$ and scaling factors of $f$.

\subsection{ACAT2P methods}	  
The expression of the  Adaptive Compact Approximate Taylor Method (ACAT$2P$) of maximal order $2P$ for a scalar conservation law is then given by:

\begin{equation}
\label{adaptivecons}
u_i^{n+1} = u_i^n + \frac{\Delta t}{\Delta x}\left( {F}^{A}_{i-1/2} - {F}^{A}_{i+1/2}\right).
\end{equation}

The numerical  fluxes ${F}^{A}_{i+1/2}$ are defined by \eqref{admissibleset}-\eqref{numfluxacat} where $F^*_{i+1/2}$ is the numerical flux of the FL-CAT$2$
\eqref{ACAT2_flux} and the smoothness indicators are given by \eqref{indicador1},  \eqref{indicadores_locales}. For $p = 2$
 \eqref{indicadores_locales} can be replaced by \eqref{indicadorlocal2mod}.

Observe that, by definition,  ${F}^{A}_{i+1/2}$  reduces to:
\begin{itemize}
	\item a first order flux if  $\psi^1_{i+1/2}=0$ and $\psi^p_{i+1/2}=0$ for all $p=2,\ldots,P;$ 
	\item a second order flux if $\psi^1_{i+1/2}=1$ and $\psi^p_{i+1/2} \approx 0$ for all $p=2,\ldots,P;$ 
	\item $2p_s$-order flux if $\psi^{p_s}_{i+1/2} \approx 1$.
\end{itemize}


Furthermore, if $p_s = P$, then  ACAT$2P$ coincides with  CAT$2P$ which has $2P$-order accuracy and is $L^2$-stable  under  CFL$\leq 1.$ 

Let us suppose that $f$ is smooth and has an isolated critical point $x^*$ of order $k$ in $S_1 = \{ x_i, x_{i+1} \}$. Then:
\begin{itemize}
    \item If $k < 2(P-1)$ the smoothness indicator $\psi^P_{i+1/2}$ is close to one and the maximum allowed stencil $S_P$ is used, so that the local accuracy of the method is $2P$.
    \item If $k > 2(P-1)$ then all the smoothness indicators fail, so that the first order robust numerical method will be used.
    Nevertheless in this case,  $f^{(j)}(x^*) = 0$ for $ j = 1, \dots, 2P - 1$ so that, when the local error of the first order method is estimated through Taylor expansions, only terms of order $O(\Delta x^{2P})$ or bigger will remain. Therefore, in this case the local accuracy of the method is again $2P$.
    \item If $k = 2(P-1)$ again the smoothness indicators fail and the first order robust numerical method will be used.
    Since in this case,  $f^{(j)}(x^*) = 0$ for $ j = 1, \dots, 2P - 2$  the local error of the first order method is  of order $2P-1$.
\end{itemize}
Summing up, the local accuracy of the method close to a critical point is always $2P$ with the only exception of critical points of order
$2P-2$: in that case, the order of accuracy will be reduced by one. This order reduction could be avoided by introducing optimal smoothness
indicators in the spirit of \cite{ZBBM2019},\cite{BBMZO3}.

\subsection{Systems of conservation laws}\label{ss:systems}
For systems of conservation laws \eqref{cons_law_0} with $m > 1$ the expression of the ACAT$2P$ method is the same as in the scalar case: the only difference is the computation of the smoothness indicators. In the case of systems, smoothness indicators are first computed for every variable:
$$\psi^{j,p}_{i+1/2}, \quad p = 1, \dots, P,$$
where
\begin{itemize}
    \item $\psi^{j,1}_{i+1/2}$ is obtained by applying the smoothness indicator \eqref{indicador1}, \eqref{local_indicators} to the $j$th component of the numerical solutions
    $\{u^{j,n}_i\}$.
    \item $\psi^{j,p}_{i+1/2}$, $p >2$  is obtained by applying the smoothness indicator \eqref{indicadores_locales} to the $j$th component of the numerical solutions
    $\{u^{j,n}_i\}$.
        \item $\psi^{j,2}_{i+1/2}$ is obtained by applying the smoothness indicator \eqref{indicadores_locales} or \eqref{indicadorlocal2mod} to the $j$th component of the numerical solutions
    $\{u^{j,n}_i\}$.
\end{itemize}

Once these scalar  smoothness indicators have been computed, we define

$$
\psi^p_{i+1/2} = \min_{j=1, \dots, m } \psi^{j,p}_{i+1/2}, 
$$
so that the selected stencil is the one of maximal length among  those in which all the variables are smooth. 

\begin{remark}
Standard WENO schemes applied componentwise usually produce oscillatory solutions near shock discontinuities. To alleviate this problem, it is possible to perform a WENO reconstruction on the characterisctic variables, as described in \cite{Qiu2002}. 
This technique reduces the oscillations but dramatically increases the computational cost.
Here we do not feel the need of such a procedure, since our reconstructions are usually much less oscillatory than componentwise WENO.
\end{remark}

\section{Two-dimensional problems} 
\begin{figure}[h!]
	\centering
	\begin{picture}(80,175)
	\put(10,4){\makebox(60,165)[c]{		
			\includegraphics[scale= 0.6]{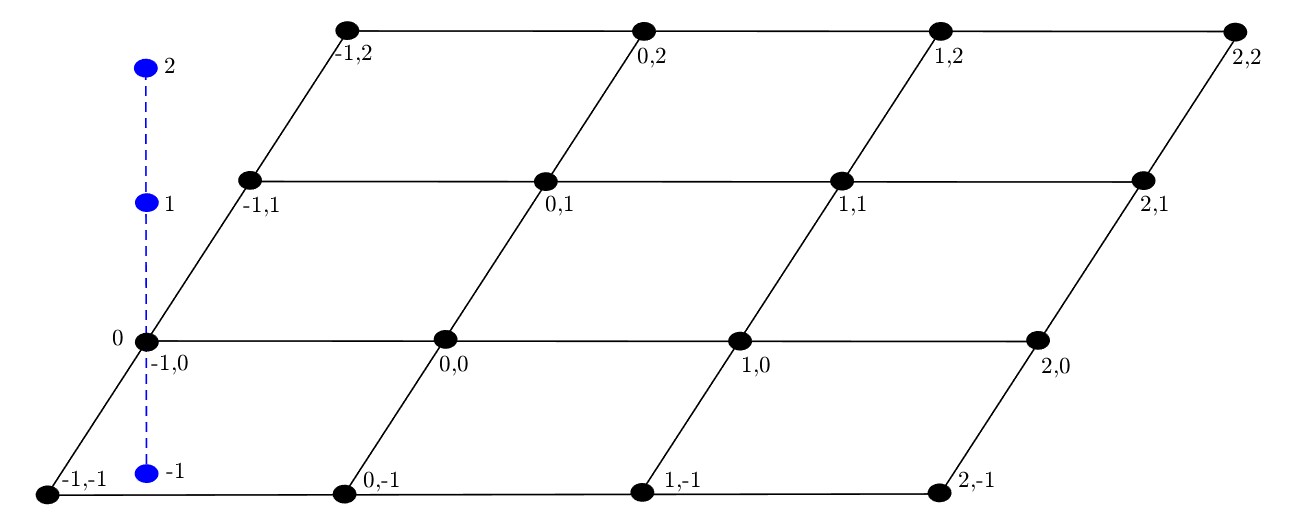}}}
	\end{picture}
	\caption{Stencil $S_2$ centered in $\mathbf{x}_\bhalf = (0.5\Delta x, 0.5\Delta y$)}
	\label{2D_CAT_1}
\end{figure}

In this section we focus on the extension of ACAT methods to non-linear two-dimensional systems of hyperbolic conservation laws
\begin{equation}
\label{2Dequ}
{u}_t + {f}({u})_x + {g}({u})_y=0.
\end{equation}
The following multi-index notation will be used:
$$
\bi = (i_1, i_2) \in \mathbb{Z} \times \mathbb{Z},
$$
and
$$
\bzero = (0,0),\quad \bone = (1,1), \quad \bhalf= (1/2, 1/2), \quad \bef = (1,0), \quad \bes = (0,1).
$$
We consider Cartesian meshes with nodes
$$\mathbf{x}_{\bi} = (i_1 \Delta x, i_2 \Delta y).$$
Using this notation, the general form of the CAT$2p$ method will be as follows:
\begin{equation}
\label{2Dscheme}
u_\bi^{n+1}=u_\bi^n + \frac{\Delta t}{\Delta x}\left[ \mathcal{F}_{\bi -  \frac{1}{2}\bef}^p-\mathcal{F}_{\bi +  \frac{1}{2}\bef}^p\right] +
\frac{\Delta t}{\Delta y}\left[ \mathcal{G}_{\bi -  \frac{1}{2}\bes}^p-\mathcal{G}_{\bi +  \frac{1}{2}\bes}^p\right],  
\end{equation}
where the numerical fluxes $\mathcal{F}_{\bi +  \frac{1}{2}\bef}^p$,  $\mathcal{G}_{\bi +  \frac{1}{2}\bes}^p$ will be computed using the values of the numerical solution $u_{\bi}^n$ in the $p^2$-point stencil centered at $\mathbf{x}_{\bi + \bhalf} = 
((i_1 + 1/2)\Delta x, (i_2 + 1/2)\Delta y)$ 
$$
S_p = \{ \mathbf{x}_{\bi + \bj}, \quad \bj \in \mathcal{I}_p  \}, 
$$
where 
$$
\mathcal{I}_p =\{ \bj = (j_1, j_2) \in \Z \times \Z, \quad -p +1 \leq j_k \leq p, \quad k = 1,2 \}.
$$
See Figure \ref{2D_CAT_1} for an example. 

For instance, the expression of the CAT2 numerical flux is as follows:
\begin{align}
\label{2Dfluxes}
F_{\bi +  \frac{1}{2}\bef}^* = & \frac{1}{4}\left( \tilde{f}^{1,n+1}_{\bi,\bzero} + \tilde{f}^{1,n+1}_{\bi, \bef} + f^{n}_{\bi} +f^{n}_{\bi + \bef}\right), \\
G_{\bi +  \frac{1}{2}\bes}^* = & \frac{1}{4}\left( \tilde{g}^{1,n+1}_{\bi,\bzero} + \tilde{g}^{1,n+1}_{\bi, \bes} + g^{n}_{\bi} +g^{n}_{\bi + \bes}\right), 
\end{align}
where 
\begin{align*}
\tilde{f}^{1,n+1}_{\bi, \bj} &= f\left(u_{\bi + \bj}^n +\Delta t \tilde{u}^{(1)}_{\bi, \bj}\right), \\
\tilde{g}^{1,n+1}_{\bi, \bj} &= g\left(u_{\bi + \bj}^n +\Delta t \tilde{u}^{(1)}_{\bi, \bj}\right),
\end{align*}
for $\bj= \bzero, \bef, \bes$. Furthermore,
\begin{align*}
\tilde{u}^{(1)}_{\bi, \bzero} &= -\frac{1}{\Delta x}\left( f_{\bi + \bef}^n - f_{\bi}^n\right) 
-\frac{1}{\Delta y}\left( g_{\bi + \bes}^n - g_{\bi}^n\right),  \\
\tilde{u}^{(1)}_{\bi, \bef} &= -\frac{1}{\Delta x}\left( f_{\bi + \bef}^n - f_{\bi}^n\right)
-\frac{1}{\Delta y}\left( g_{\bi + \bone}^n - g_{\bi + \bef}^n\right), \\ 
\tilde{u}^{(1)}_{\bi, \bes} &= -\frac{1}{\Delta x}\left( f_{\bi + \bone}^n - f_{\bi + \bes}^n\right) 
-\frac{1}{\Delta y}\left( g_{\bi + \bes}^n - g_{\bi}^n\right), \\
\end{align*} 
where 
$$ f_\bj^n = f(u_{\bj}^n), \quad g_\bj^n =g(u_{\bj}^n), \quad \forall \bj.$$
Observe  that $\tilde{u}^{(1)}_{\bi, \bzero} \neq \tilde{u}^{(1)}_{\bi, \bef}$ and $\tilde{u}^{(1)}_{\bi, \bzero} \neq \tilde{u}^{(1)}_{\bi, \bes}$ as opposed to the 1D case where $\tilde{u}^{(1)}_{i,0 }= \tilde{u}^{(1)}_{i, 1} $: compare with \eqref{CAT2_flux}-\eqref{CAT2_flux_2}.
The following algorithm will be used to compute the numerical fluxes of the CAT2$p$ method:

\begin{enumerate}

\item  {Define}
$$
\tilde{f}^{(0)}_{\bi,\bj}=f^n_{\bi + \bj}, \quad  \tilde{g}^{(0)}_{\bi,\bj}=g^n_{\bi + \bj}, \quad \bj \in \mathcal{I}_p.
$$

\item  {For $k = 2 \dots 2p$:}

\begin{enumerate}

\item Compute
\begin{equation*}
 \tilde{u}^{(k-1)}_{\bi,\bj} = - A^{1,j_1}_{p,0}(\tilde{f}^{(k-2)}_{\bi,(\bcdot, j_2)}, \Delta x) 
 - A^{1,j_2}_{p,0}(\tilde{g}^{(k-2)}_{\bi,(j_1, \bcdot)}, \Delta y), \quad \bj \in \mathcal{I}_p.
\end{equation*}

\item Compute 
$$
\tilde{f}^{k-1,n+r}_{\bi,\bj} = f \left(  u^n_{\bi + \bj} + 
\sum_{l=1}^{k-1} \frac{(r \Delta t)^l}{l!} \tilde{u}^{(l)}_{\bi,\bj} \right), \quad  \bj\in \mathcal{I}_p,\> r = -p+1, \dots, p.
$$

\item Compute
$$
\tilde{f}^{(k-1)}_{\bi,\bj} =   A^{k-1,0}_{p,n}( \tilde{f}^{k-1, \bcdot}_{\bi,\bj}, \Delta t),\quad  \bj \in \mathcal{I}_p.
$$

\end{enumerate}

\item Compute
\begin{eqnarray}\label{cat22Dx}
F^p_{\bi+\frac{1}{2}\bef}  &= &\sum_{k=1}^{2p} 
\frac{\Delta t^{k-1}}{k!}A^{0, 1/2}_{p,0}(\tilde{f}_{\bi,(\bcdot, 0)}^{(k-1)}, \Delta x), 
\\\label{cat22Dy}
G^p_{\bi+\frac{1}{2}\bes}  &= &\sum_{k=1}^{2p} 
\frac{\Delta t^{k-1}}{k!}A^{0, 1/2}_{p,0}(\tilde{g}_{\bi, (0, \bcdot)}^{(k-1)}, \Delta y).
\end{eqnarray}

\end{enumerate}
The notation used for the approximation of the spacial partial derivatives is the following:
\begin{eqnarray*}
 A^{k,q}_{p,j_1}(f_{\bi, (\bcdot, j_2)}, \Delta x) & =&  \frac{1}{\Delta x^k} \sum_{l = -p + 1}^p \gamma^{k,q}_{p,l} f_{\bi, (l,j_2)}\\
 A^{k,q}_{p,j_2}(g_{\bi, (j_1,\bcdot{} )}, \Delta y) & =&  \frac{1}{\Delta y^k} \sum_{l = -p + 1}^p \gamma^{k,q}_{p,l} g_{\bi, (j_1,l)}
 \end{eqnarray*}
 
 \begin{remark}
 In the last step of the algorithm above the set $\mathcal{I}_p$ can be replaced by its $(2p -1)$-point subset
 $$
 \mathcal{I}^0_p = \{ \bj = (j_1, j_2) \text{ s.t } j_1 = 0 \text{ or } j_2 = 0 \}
 $$
 since only the corresponding values of $\tilde{f}^{(k-1)}_{\bi,\bj}$ are  used to compute the numerical fluxes \eqref{cat22Dx} and
 \eqref{cat22Dy}.
 \end{remark}

Once the numerical flux of the CAT$2p$ method has been introduced, the numerical flux of ACAT2 is extended to two-dimensional problems as follows:
\begin{eqnarray}\label{ACAT2_flux_F} 
\mathcal{F}_{\bi +  \frac{1}{2}\bef}^1 &= &\psi^1_{\bi + \frac{1}{2}\bef} \, F^*_{\bi +  \frac{1}{2}\bef} +(1-\psi^1_{\bi + \frac{1}{2}\bef}) \, F^{lo}_{\bi +  \frac{1}{2}\bef},  
\\\label{ACAT2_flux_G} 
\mathcal{G}_{\bi +  \frac{1}{2}\bes}^1 &= &\psi^1_{\bi + \frac{1}{2}\bes} \, G^*_{\bi +  \frac{1}{2}\bes} +(1-\psi^1_{\bi + \frac{1}{2}\bes}) \, G^{lo}_{\bi +  \frac{1}{2}\bes},  
\end{eqnarray}
where, $F^{lo}_{\bi +  \frac{1}{2}\bef}$ and $G^{lo}_{\bi +  \frac{1}{2}\bes}$ are some robust first order methods; 
$\psi^1_{\bi + \frac{1}{2}\bef}$ and $\psi^1_{\bi + \frac{1}{2}\bes}$ are the flux limiters computed dimension by dimension. 

Finally, the expression of the ACAT$2P$ method for two-dimensional problems is
\begin{equation}
\label{ACAT2P2D}
u_\bi^{n+1}=u_\bi^n + \frac{\Delta t}{\Delta x}\left( \mathcal{F}_{\bi -  \frac{1}{2}\bef}^{A_1}-\mathcal{F}_{\bi +  \frac{1}{2}\bef}^{A_1}\right) +
\frac{\Delta t}{\Delta y}\left( \mathcal{G}_{\bi -  \frac{1}{2}\bes}^{A_2}-\mathcal{G}_{\bi +  \frac{1}{2}\bes}^{A_2}\right),  
\end{equation}
where the numerical fluxes are defined as follows: first define the set

\begin{eqnarray}\label{admissibleset_2Dx}
\mathcal{A}_1 & = & \{  p \in \{2, \dots, P \} \text{ s.t. } 
\psi^p_{\bi+\frac{1}{2}\bef \cong 1} \},\\
\label{admissibleset_2Dy}
\mathcal{A}_2 & = & \{  p \in \{2, \dots, P \} \text{ s.t. }  \psi^p_{\bi+\frac{1}{2}\bes}\cong 1\},\\
\end{eqnarray}
where $\psi^p_{\bi+\frac{1}{2}\bef}$, $\psi^p_{\bi+\frac{1}{2}\bes}$ are the smoothness indicators introduced in Section \ref{ss:smoothness} computed dimension by dimension. Then define:
\begin{eqnarray}\label{numfluxacat2dx}
F^{A_1}_{\bi+\frac{1}{2}\bef} &=& \begin{cases}
F^{*}_{\bi+\frac{1}{2}\bef} & \text{if} \> \mathcal{A}_1 = \emptyset;\\
F^{p_1}_{\bi+\frac{1}{2}\bef} & \text{where} \,\, p_1 = \max(\mathcal{A}_1) \text{ otherwise;} \\
\end{cases}
\\\label{numfluxacat2dy}
G^{A_2}_{\bi+\frac{1}{2}\bes} &=& \begin{cases}
G^{*}_{\bi+\frac{1}{2}\bes} & \text{if} \,\,\mathcal{A}_2 = \emptyset ;\\
G^{p_2}_{\bi+\frac{1}{2}\bes} & \text{where} \,\,p_2 = \max(\mathcal{A}_2) \text{ otherwise.} \\
\end{cases}
\end{eqnarray}

Observe that, since the smoothness indicators are computed dimension by dimension, a rectangular stencil
$$
S_{p_1, p_1} = \{ \mathbf{x}_{\bi, \bj}, \quad i_1-p_1 + 1 \leq j_1 \leq i_1 + p_1, \quad i_2-p_2 + 1 \leq j_2 \leq i_2 + p_2 \},
$$
is used in practice to compute the numerical fluxes $ F^{p_1}_{\bi+\frac{1}{2}\bef}$, $G^{p_2}_{\bi+\frac{1}{2}\bes}$. The extension of CAT methods to such rectangular stencils is straightforward.

\section{Numerical experiments}
In this section we apply ACAT$2P$ methods to several 1D and 2D problems: the 1D linear transport equation, Burgers equation, and the 1D and 2D Euler equation for gas dynamic.
The Super Bee flux limiter \cite{Roesuperbee} is used in FL-CAT2  and the smoothness indicators \eqref{indicadores_locales} are used for $p \geq 2$: no loss of precision for first order critical points has been observed in any of the test problems considered here due to the use of $\psi^2_{i+1/2}$. Fornberg's algorithm \cite{Fornberg1} is used to compute the coefficients of the numerical differentiation formulas. ACAT methods will be compared with the Lax-Friedrichs (LF), HLL first order schemes and with
WENO($2p+1$) finite difference methods based on the Lax-Friedrichs splitting (see \cite{Shu1997}) combined with SSPRK3 (\cite{Gottlieb2011}) for the time discretization.
The order and the number of points of their stencils in 1d are recalled in Table \ref{tmethods}. Since ACAT$2P$ reduces to CAT$2P$ and the order of accuracy of the latter have been checked in \cite{CP2019}, no test order will be considered here: interested readers are referred to that work.

\begin{table}[h]
\centering
\begin{tabular}{|c||c||c|}\hline
Method & Stencil & Order \\ \hline \hline 
LF &  3  & 1  \\ \hline
HLL & 3 & 1  \\ \hline
ACAT2 or FL-CAT2 & 3 & 2  \\ \hline
ACAT$2P$ & $2P+1$ & $2P$  \\ \hline
WENO($2p+1$)-RK3 & $2p+1$ & $2p+1$\\ \hline
\end{tabular}
\caption{Numerical methods: order of accuracy and number of points of the stencils for 1d problems.\label{tmethods}}
\end{table}

\subsection{1D linear transport equation}

Let us consider the linear scalar conservation law


\begin{equation}\label{cons_scalar}
u_t+ u_x=0.
\end{equation}
with initial condition:

\begin{equation}\label{condsuave_1}
u_0(x)=\frac{1}{2}\sin(\pi x).
\end{equation} 

\begin{figure}[!ht]
	\centering	
	\includegraphics[width=\textwidth,height=8 cm]{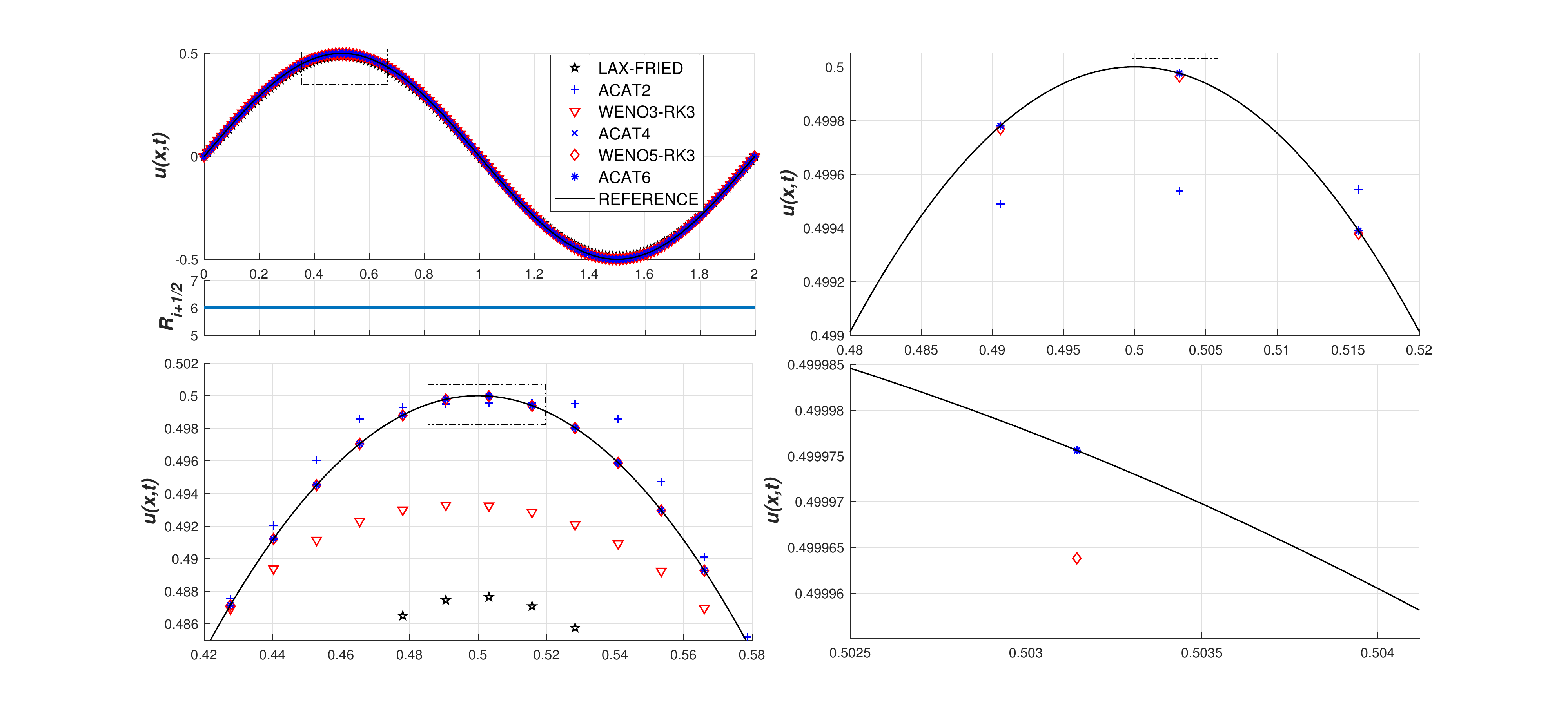}
	\vspace{-0.8 cm}
	\caption{Transport equation with initial condition (\ref{condsuave_1}).  Numerical solution at $t = 4$: general view (\textit{left-up}); order of accuracy for ACAT$6$ (\textit{sub-frame}); consecutive zooms close to the local maximum ( \textit{left-down}, \textit{right-up} and \textit{right-down}).}
	\label{test1}
\end{figure}

\begin{figure}[!ht]
	\centering	
	\includegraphics[width=\textwidth,height=8 cm]{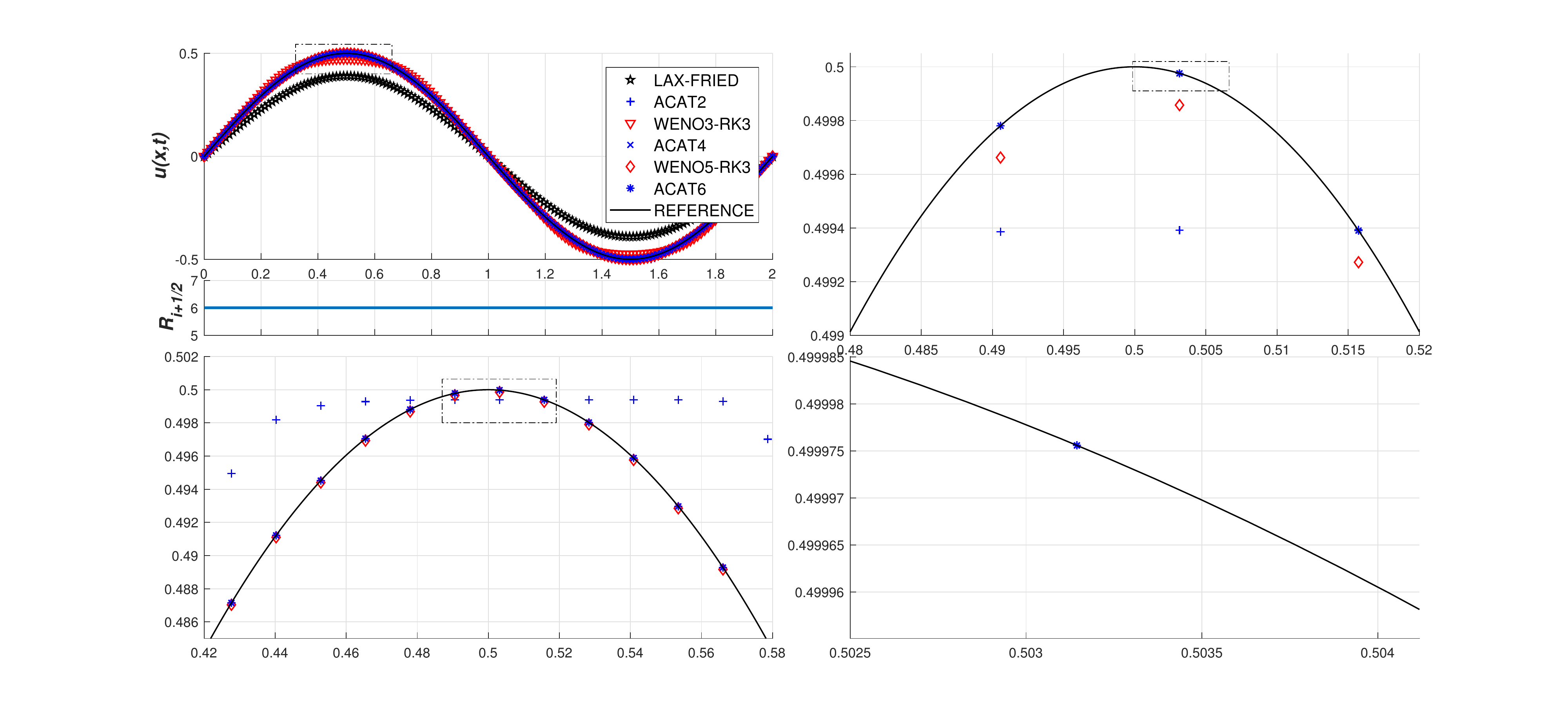}
	\vspace{-0.8 cm}
	\caption{Transport equation with initial condition (\ref{condsuave_1}).  Numerical solution at $t = 40$: general view (\textit{left-up}); local order of accuracy for ACAT$6$ (\textit{sub-frame});consecutive zooms close to the local maximum ( \textit{left-down}, \textit{right-up} and \textit{right-down}).}
	\label{test1_2}
\end{figure}

We solve numerically this problem in the spatial interval $[0,2]$, using a $160$-mesh points, CFL$=0.9$, and periodic boundary conditions.

Figure \ref{test1}  and \ref{test1_2} show the numerical solutions at time $t = 4$ and $t = 40$ respectively. Zooms of an interest area are included, in which the loss of accuracy with time for the  lower order methods can be clearly seen. As it can be observed, the numerical solutions of ACAT$4$ and ACAT$6$ match the exact solution at both times while ACAT$2$ is more diffusive near the critical points. This loss of accuracy close to the critical points can also be observed for WENO-RK methods, although this drawback can be overcome by using optimal weights in the WENO reconstructions: see  \cite{ZBBM2019},\cite{BBMZO3}

\begin{figure}[!ht]
	\small
	\setlength{\unitlength}{1mm}
	\centering	
	\includegraphics[width=\textwidth,height=8 cm]{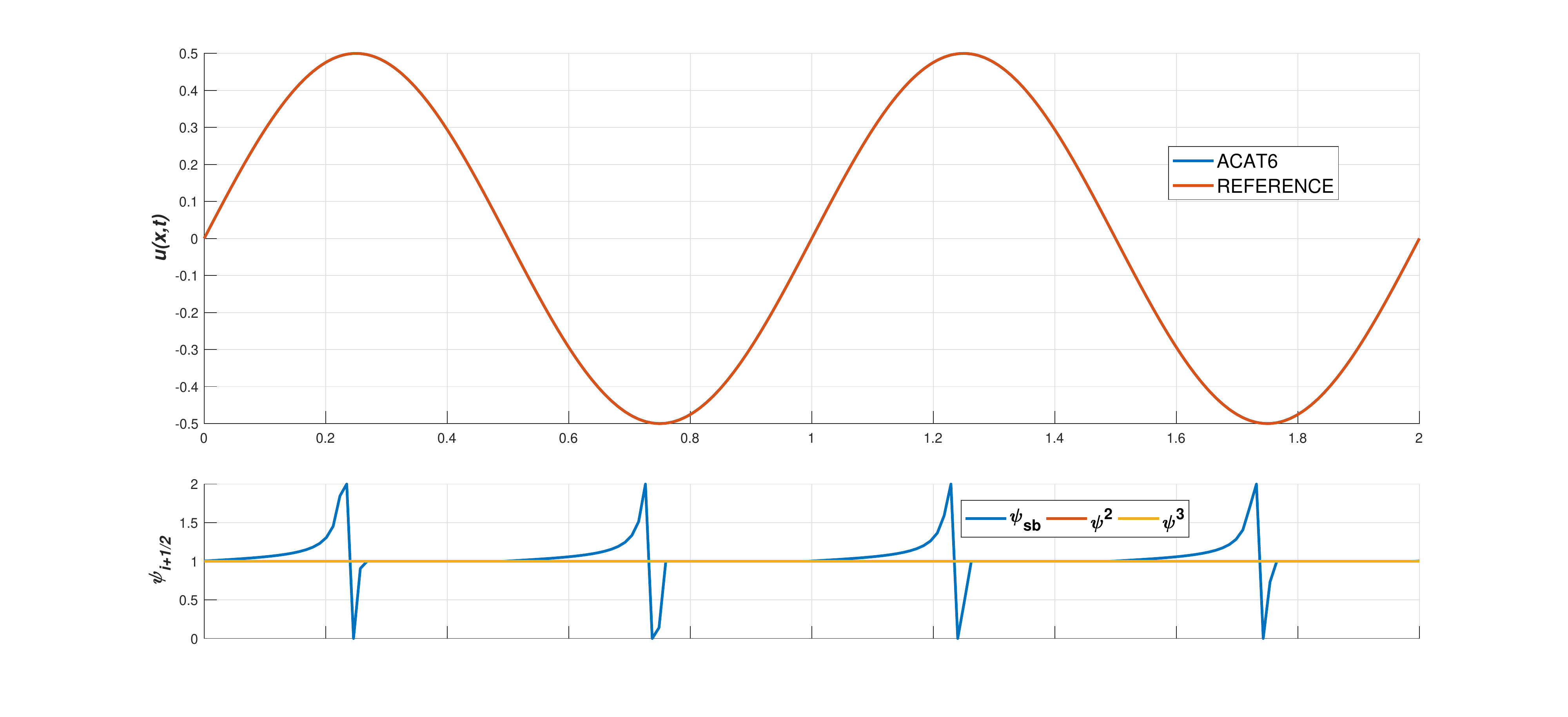}
	\vspace{-0.4 cm}
	\caption{Transport equation with initial condition (\ref{condsuave_1bis}). 
	Solution obtained with ACAT6 at time 4 (\textit{up}) and graphs of the smoothness indicators  $\psi_{sb}$, $\psi^2$ and $\psi^3$ (\textit{down}).}
	\label{test2}
\end{figure}

The loss of accuracy of ACAT$2$ close to the critical points compared to ACAT$4$ or $6$  is due to the fact that, while the smoothness indicators 
$\psi^2_{i+1/2}$ and $\psi^3_{i+1/2}$ are always close to one, the Superbee flux limiter
$\psi_{sb, i+1/2}$ detects a discontinuity at the critical points and the first order methods is then locally used: to make this clear, Figure \ref{test2} (up) shows the solution obtained with  ACAT6 at time $t = 4$ for (\ref{cons_scalar}) with initial condition 

\begin{equation}\label{condsuave_1bis}
u_0(x)=\frac{1}{2}\sin(2\pi x)
\end{equation} 
in the interval $[0,2]$ using again a 160-point mesh, CFL = 0.9, and periodic boundary conditions. Figure \ref{test2} (down) shows the graph of the three smoothness indicators. 

We consider next equation (\ref{cons_scalar}) with a piecewise continuous initial condition 
\begin{equation} \label{square_step_test}
u_0(x) = \begin{cases}
1 \quad \;\;\;\mathrm{if}\quad \frac{1}{2}\le x\le 1;\\
0 \quad\;\; \;\mathrm{if}\quad 0\le x < \frac{1}{2} \quad \mathrm{or} \quad \frac{3}{2}<x\le 2;\\
-1 \quad \mathrm{if} \quad 1<x\le\frac{3}{2}. 
\end{cases}
\end{equation} 
We solve numerically this problem in the spatial interval $[0,2]$, using again a $160$-mesh points, CFL=0.9, and periodic boundary conditions. 
\begin{figure}[!ht]
	\small
	\setlength{\unitlength}{1mm}
	\centering	
	\includegraphics[width=\textwidth,height=9 cm]{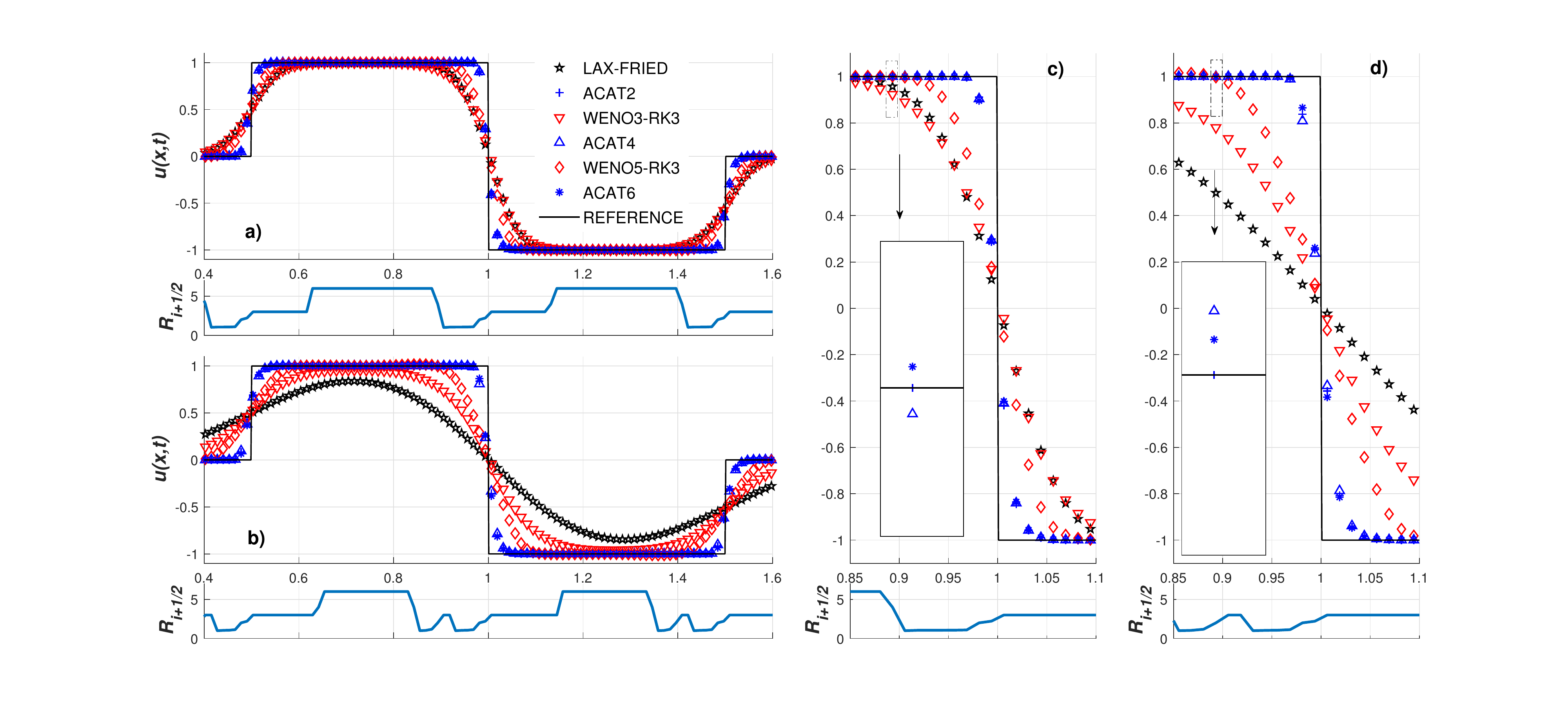}
	\vspace{-0.8 cm}
	\caption{Transport equation with initial condition (\ref{square_step_test}).  Numerical solutions at $t = 2$ (\textit{a}) and at  $t = 20$ (\textit{b}). Zooms of the numerical solutions at time $t = 2$ (\textit{c}) and $t = 20$ (\textit{d}).  Sub-frames: local order of accuracy for ACAT6.}
	\label{test3}
\end{figure}

Figure \ref{test3} shows solutions from ACAT$2P$, $P=2,4,6$ and WENO$q$-RK3, $q = 3,5$  after 2 and 20 seconds. As it can be observed, ACAT methods capture better the discontinuity than WENO-RK schemes. In this case, ACAT$4$ and ACAT$6$  reduce to ACAT$2$ at the discontinuities due to the order adaption technique.  WENO methods give accurate solutions for short times  but spurious oscillations appear with time due to the choice CFL = 0.9.

\subsection{Burgers equation}
Let us consider the Burgers equation 
\begin{equation}\label{burgers}
u_t+ \left(\frac{u^2}{2}\right)_x=0,
\end{equation}
with initial condition (\ref{condsuave_1}).
The problem is numerically solved in the interval $[0,2]$ using an uniform mesh with $160$, CFL$=0.9$, and periodic boundary conditions. A reference solution has been computed with the  Lax-Friedrichs method using $1400$-point mesh. 
\begin{figure}[!ht]
	\setlength{\unitlength}{1mm}
	\centering	
	\includegraphics[width=\textwidth,height=8 cm]{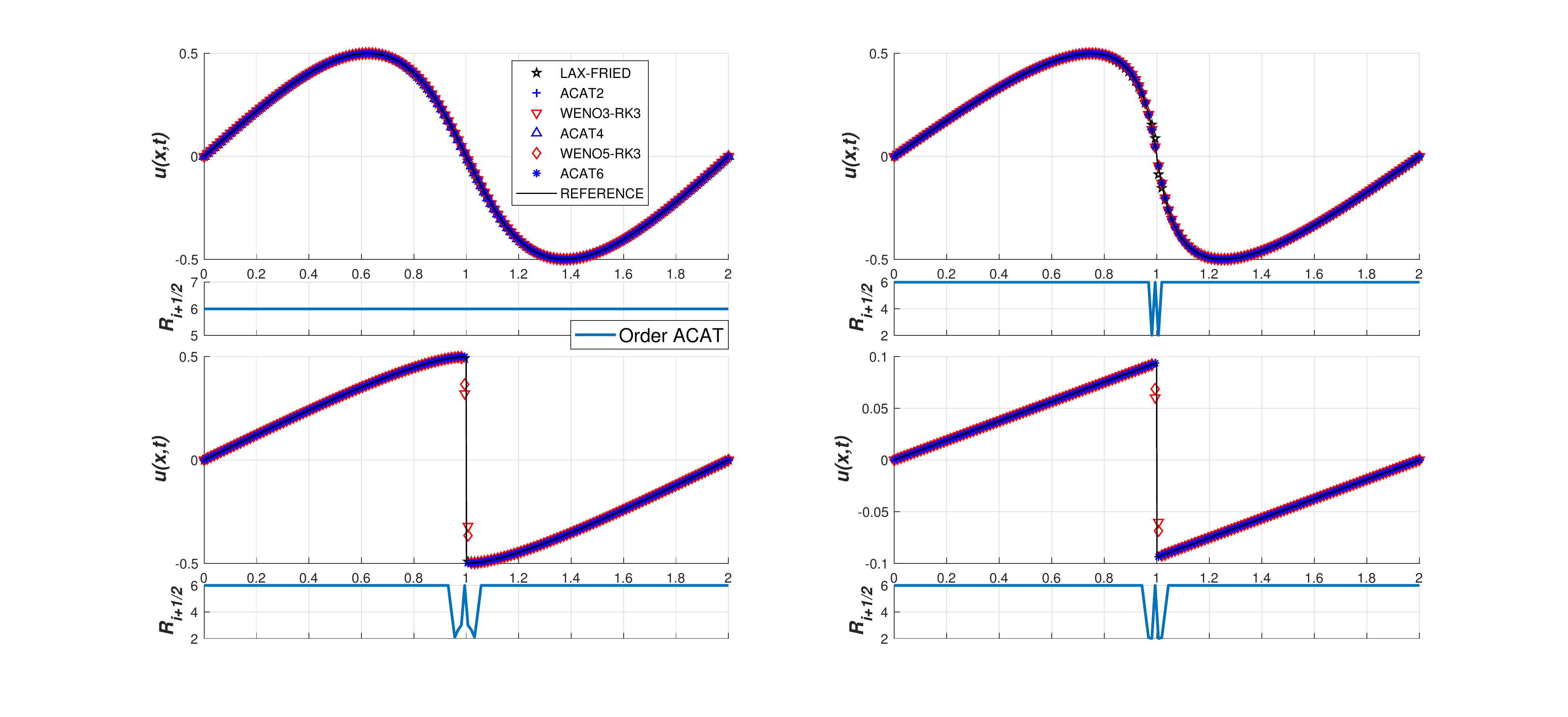}
	\vspace{-0.9 cm}
	\caption{ Burgers equation with initial condition (\ref{condsuave_1}). 
	Numerical solutions 
	obtained at times  $t = 0.25$ (\textit{left-up}), 
	$t = 0.5$ (\textit{right-up}), $t = 1$ (\textit{left-down}), and  $t = 10$ (\textit{right-down}).  Sub-frames: local accuracy order for ACAT6.}
	\label{test41}
\end{figure}

\begin{figure}[!ht]
	\setlength{\unitlength}{1mm}
	\centering	
	\includegraphics[width=\textwidth,height=8 cm]{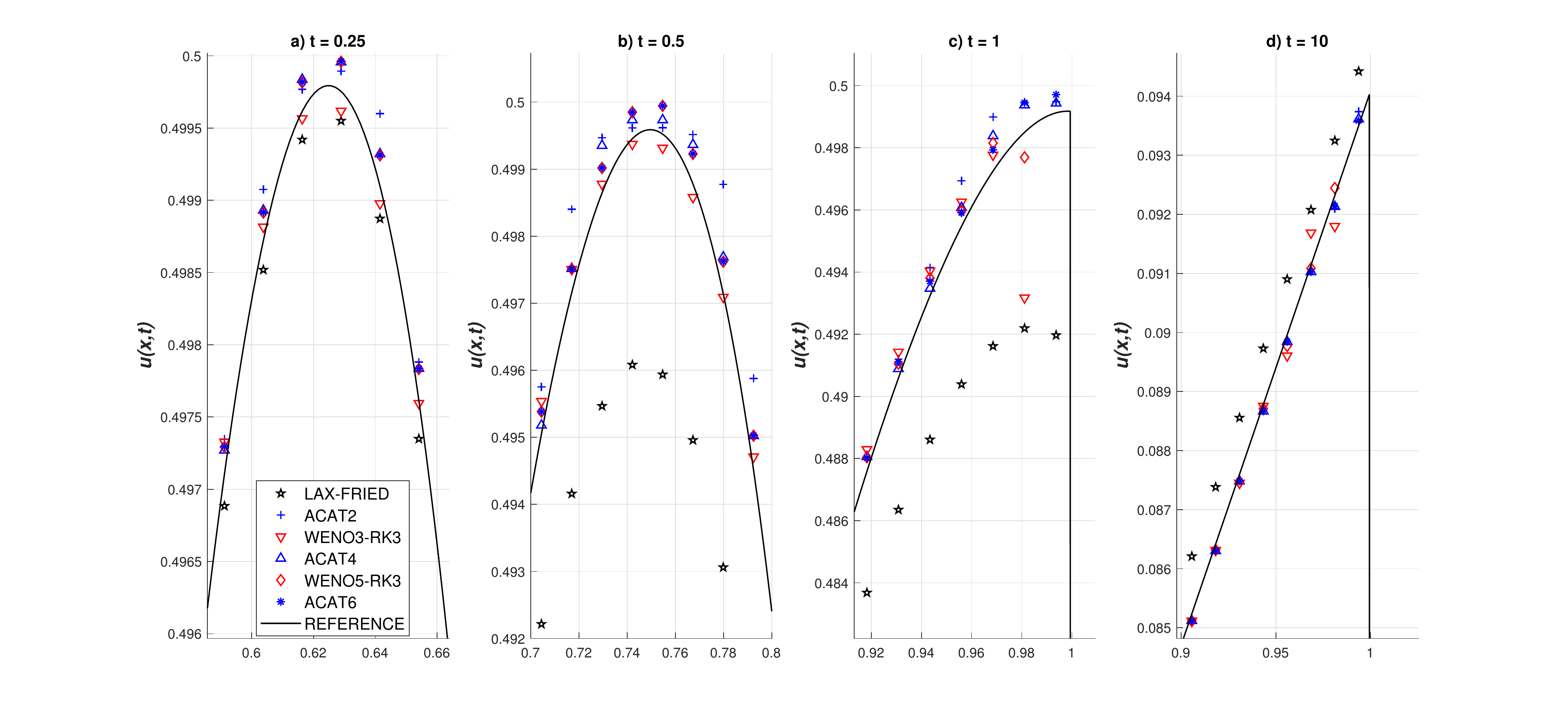}
	\vspace{-0.9 cm}
	\caption{ Burgers equation with initial condition (\ref{condsuave_1}). 
	Zoom of the numerical solutions obtained at times $t = 0.25$ (\textit{a}), $t = 0.5$ (\textit{b}), $t = 1$ (\textit{c}), and  $t = 10$ (\textit{d}). Sub-frames: local order of accuracy for ACAT6.}
	\label{test42}
\end{figure}

Figures  \ref{test41} and  \ref{test42} show respectively the general view and a zoom of the numerical solutions obtained with the different methods at times $t=\{0.25, 0.5, 1, 10\}$. 
The local order of accuracy of ACAT6 is also shown: as it can be seen, this method reduces to the first order one only at the shock once it has been generated. 

\subsection{1D Euler equations} 
Let us now consier the $1D$ Euler equations for gas dynamics
\begin{align}\label{1deuler} 
&u_t +  f(u)_x = 0,
\end{align}
with
\begin{align}
u = \left[  \begin{array}{c}
\rho \\
\rho v \\
E\\
\end{array} \right] ,\quad f(u) = \left[  \begin{array}{c}
\rho v \\
p + \rho v^2 \\
v(E+p)\\
\end{array} \right], 
\end{align}
where $\rho$ is the density measured in $Kg/m^3$; $v$, the velocity in $m/s$; $E$ the total energy per unit volume in $Kg/(ms^2)$; and $p$ is the pressure in Pascal $Pa.$ We assume an ideal gas with the  equation of state
\begin{equation}
p(\rho, e) = (\gamma - 1)\rho e,
\end{equation}
being $\gamma$ the ratio of specific heat capacities of the gas taken as 1.4 and $e$ is the internal energy per unit mass is related to $E$ by:
\begin{equation}
E = \rho (e+0.5 v^2).
\end{equation}
We consider three Riemann problems for (\ref{1deuler}):  the Sod problem \cite{Sod1978}, the Einfeldt problem  \cite{Einfeldt1991}, and the right blast wave  Woodward and Colella problem \cite{Wood1984}. In all the cases: the initial discontinuity is placed at $x = 0.5$,  the equations are numerically solved at the spatial interval $[0,1]$ and the exact solution is provided by the HE-E1RPEXACT solver introduced in \cite{Toro2009book}. The CFL parameter is set to 0.8 and outflow-inflow boundary conditions are considered. 

\begin{itemize}
	
	\item{ The Sod problem:} the initial condition is
	
	\begin{equation} \label{Sodincond}
	(\rho, v, p)
	= \left\{
	\begin{array}{ll}
	\displaystyle (1,0,1) & \mbox {if } x < 1/2, \\
	\displaystyle (0.125,0,0.1) & \mbox {if } x > 1/2.
	\end{array}\right.
	\end{equation}
	\begin{figure}[!ht]
		\setlength{\unitlength}{1mm}
		\centering	
		\includegraphics[width=\textwidth,height=8 cm]{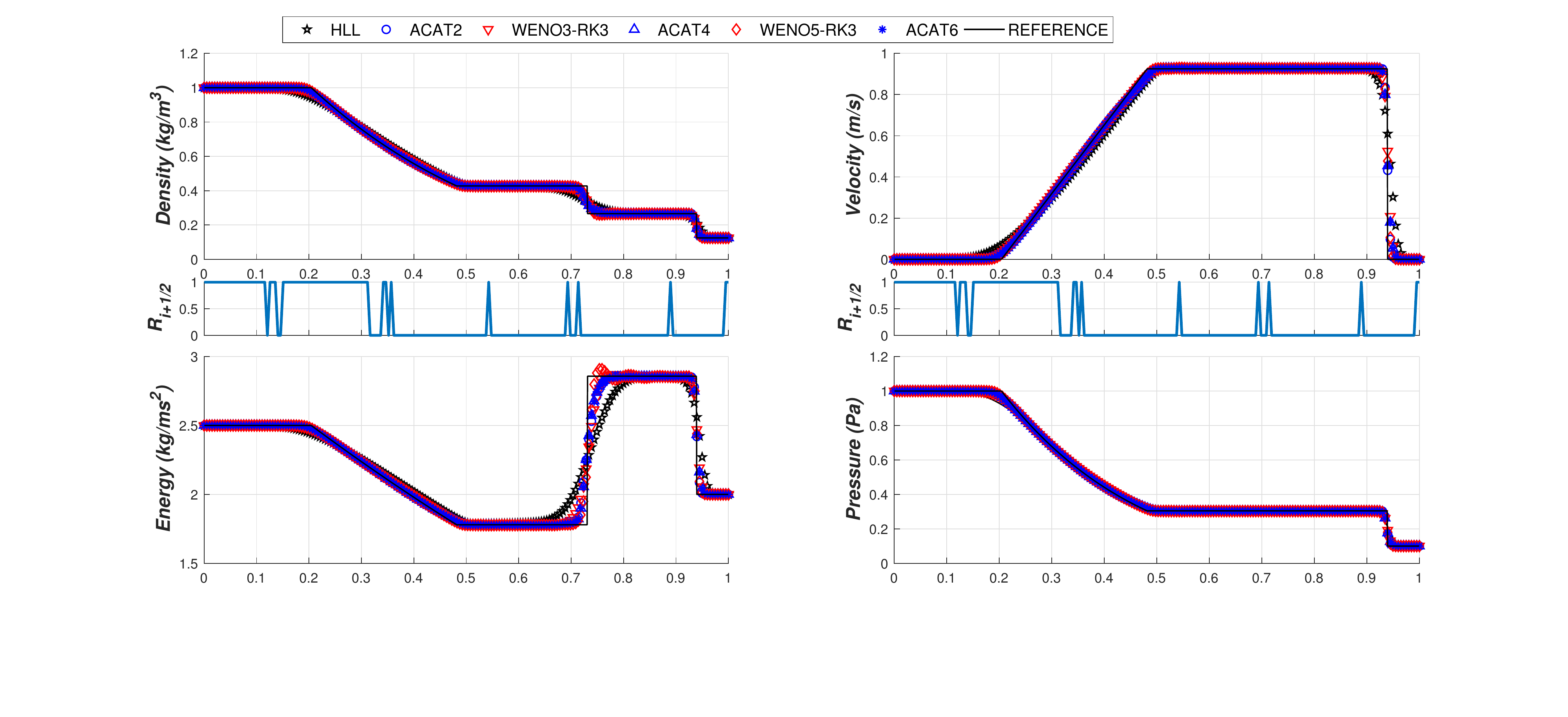}
		\vspace{-1.1 cm}
		\caption{1D Euler equations: the Sod problem.  Numerical solutions at $t = 0.25$ using CFL$=0.8$ and $200$ points: density (\textit{left-up}), velocity  (\textit{right-up}), internal energy (\textit{left-down}),  pressure (\textit{right-down}).  Sub-frames: local order of accuracy for ACAT6.}
		\label{test51}
	\end{figure}
	
	The solution involves a rarefaction wave, a contact discontinuity and a shock. We compare  the numerical solutions with the exact one: see \cite{Toro2009book}. 
	\begin{figure}[!ht]
		\setlength{\unitlength}{1mm}
		\centering	
		\includegraphics[width=\textwidth,height=8 cm]{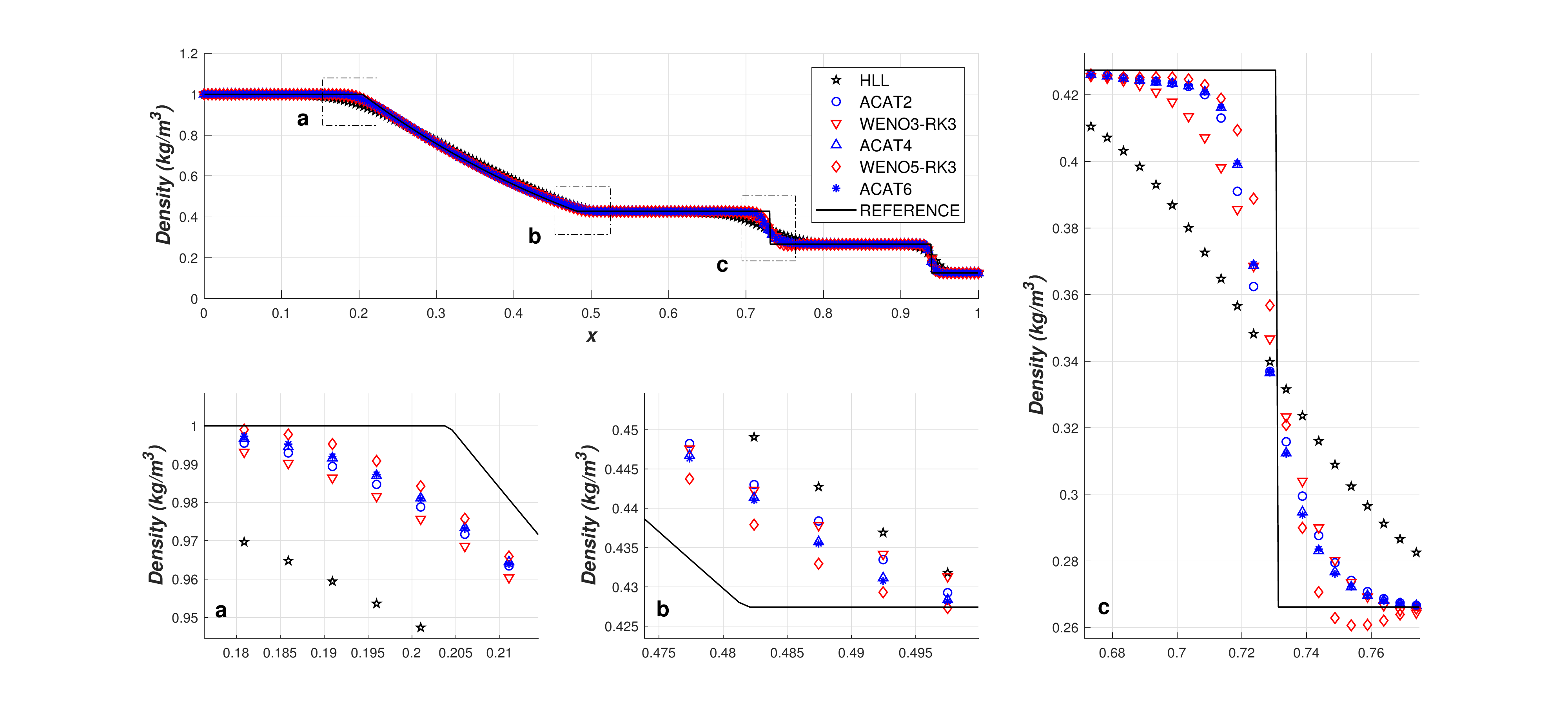}
		\vspace{-1.1 cm}
		\caption{1D Euler equations: the Sod problem.  Numerical density at $t = 0.25$ using CFL$=0.8$ and $200$ points: general view and zooms close to the points \textit{a},\textit{b}, \textit{c} and \textit{d}.}
		\label{test52}
	\end{figure}
	\begin{figure}[!ht]
		\setlength{\unitlength}{1mm}
		\centering	
		\includegraphics[width=\textwidth,height=8 cm]{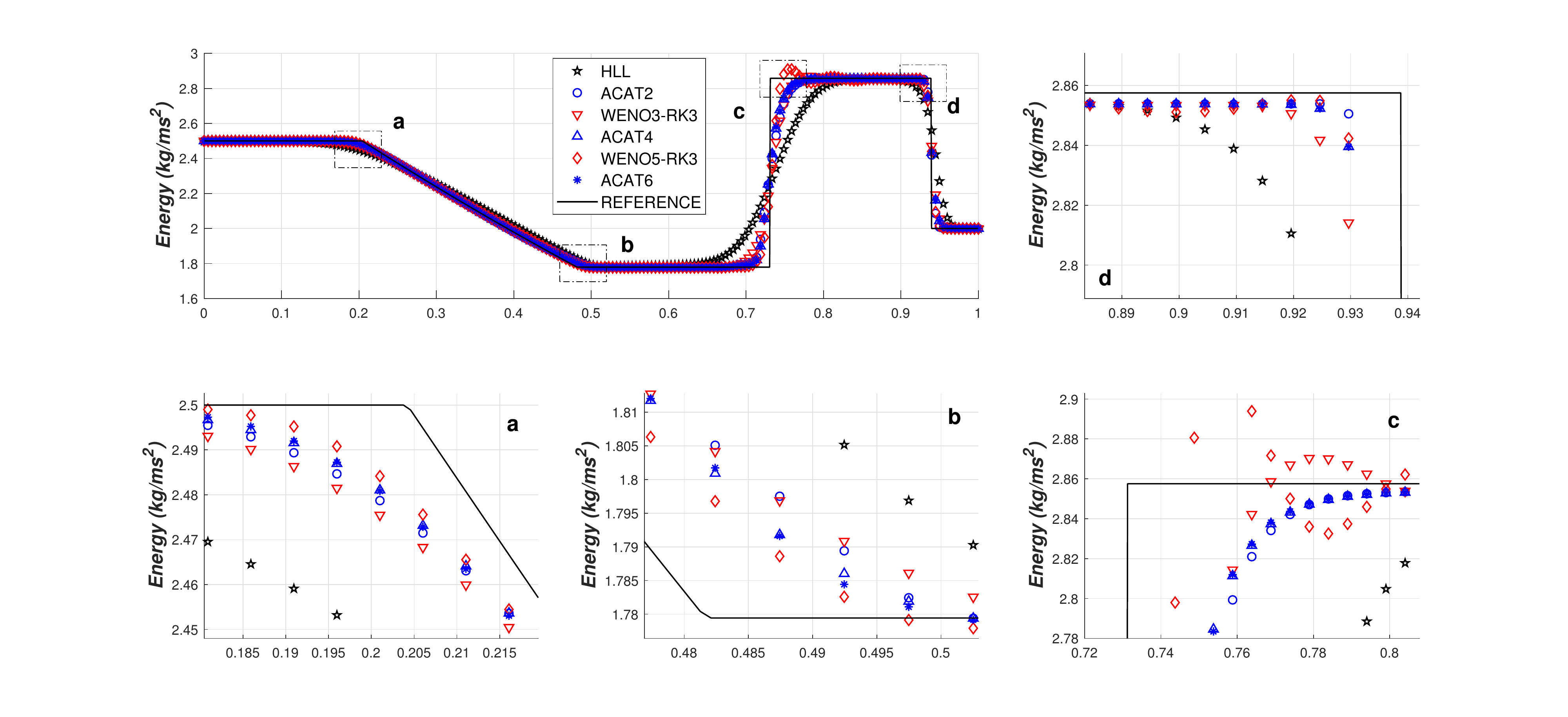}
		\vspace{-1.1cm}
		\caption{1D Euler equations: the Sod problem.  Numerical internal energy at $t = 0.25$ using CFL$=0.8$ and $200$ points: general view and zooms close to the points \textit{a},\textit{b}, \textit{c} and \textit{d}.}
		\label{test53}
	\end{figure}
	
	Figure \ref{test51}  shows the solutions provided by ACAT2-4-6 and WENO3-5 for density, velocity, internal energy and pressure $p$, using a 200-point mesh. The local accuracy of ACAT$6$ is also shown. 
	Zooms of the behaviour of the numerical densities can be observed in Figure \ref{test52}. As it can be seen in zooms \textit{a} and \textit{b}, WENO5-RK3 gives sharper but more oscillatory solutions than ACAT methods. Moreover, increasing the accuracy order for ACAT methods we obtain sharper results. Similar conclusions for the internal energy can be drawn: see Figure \ref{test53}. 
	
	\item{123 Einfeldt problem:}	the initial condition is
	\begin{equation} \label{t2incond}
	(\rho, v, p)
	= \left\{
	\begin{array}{ll}
	\displaystyle ( 1.0,-2.0,0.4) & \mbox {if } x < 1/2, \\
	\displaystyle ( 1.0, 2.0,0.4) & \mbox {if } x > 1/2.
	\end{array}\right.
	\end{equation}
		
		\begin{figure}[!ht]
		\setlength{\unitlength}{1mm}
		\centering	
		\includegraphics[width=\textwidth,height=8 cm]{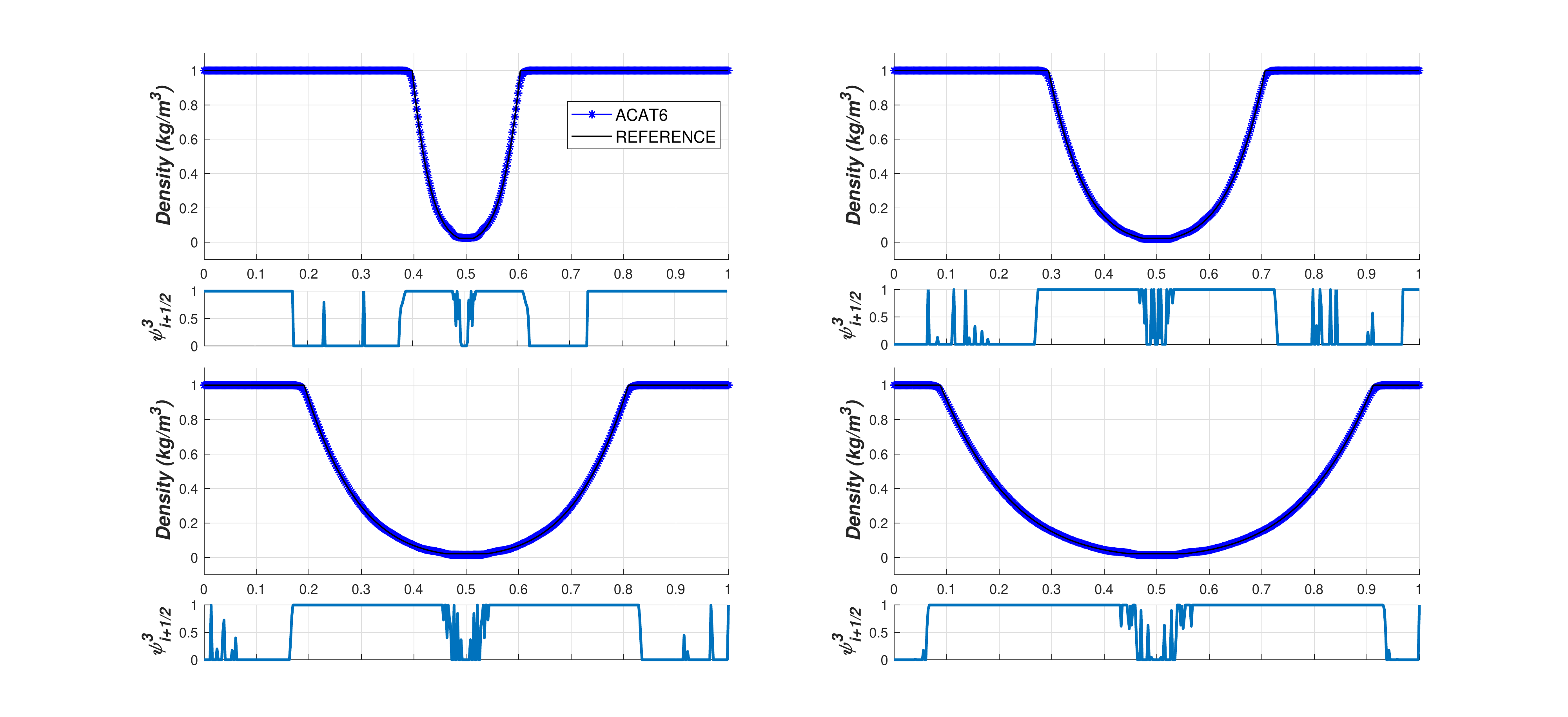}
		\vspace{-1.1 cm}
		\caption{ 1D Euler equations: the 123 Einfeldt problem. Numerical solutions at $t_s = 0.15$ using CFL$=0.8$ and $200$ points. Density obtained with ACAT6 and graph of the  smoothness indicator $\psi^3$ 
		for $t = t_s/4$ (\textit{left-up}), $t_s/2$ (\textit{right-up}), $3t_s/4$ (\textit{left-down}), $t_s$ (\textit{right-down}), with $t_s = 0.15$.}
		\label{test61}
	\end{figure} \vspace{-0.5cm}
	
	\begin{figure}[!ht]
		\setlength{\unitlength}{1mm}
		\centering	
		\includegraphics[width=\textwidth,height=8 cm]{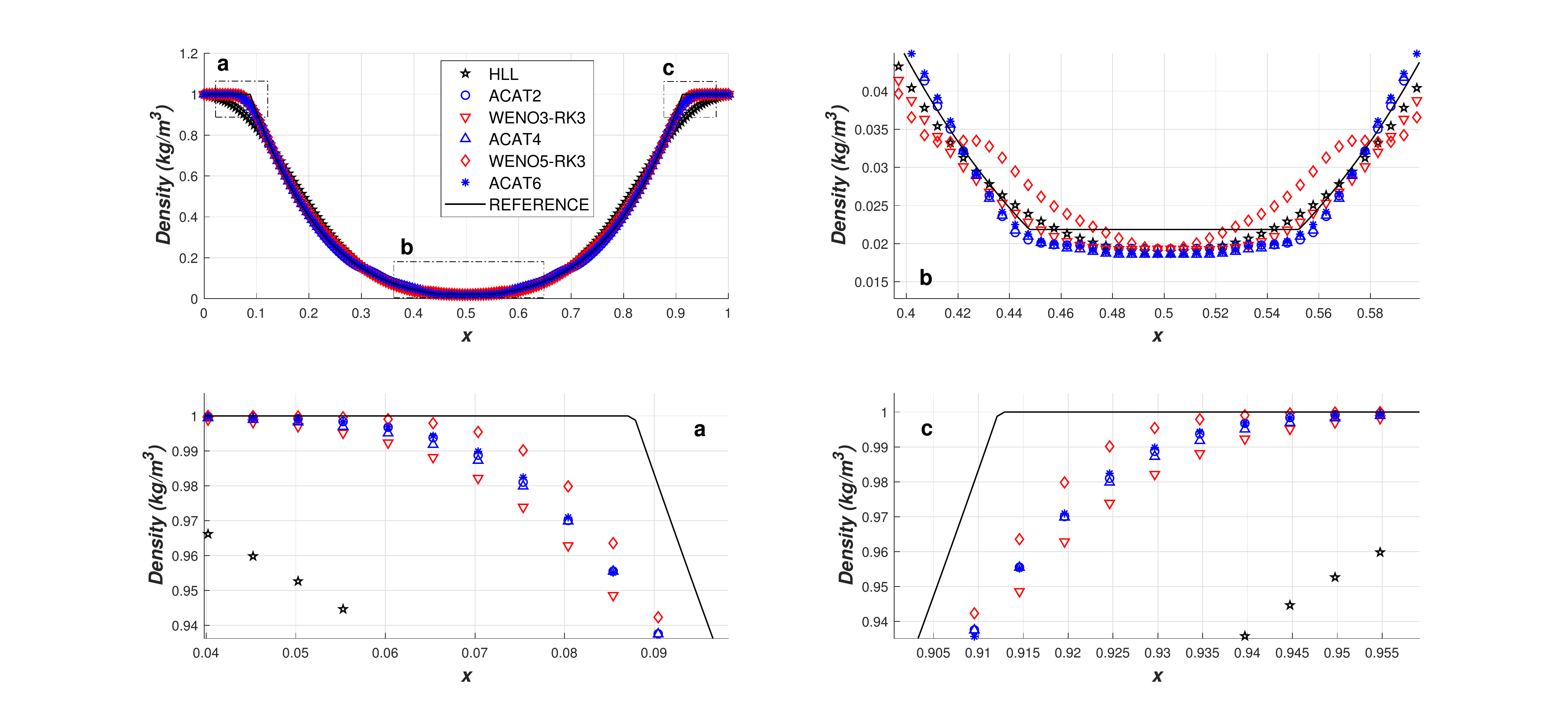}
		\vspace{-1.1 cm}
		\caption{ 1D Euler equations: the 123 Einfeldt problem.  Numerical solutions at $t = 0.15$ using CFL$=0.8$ and $200$ points: general view (\textit{left-up}) and zooms close to the points \textit{a} (\textit{left-down}), \textit{b}(\textit{right-up}), and  \textit{c} (\textit{right-down}). }
		\label{test62}
	\end{figure}

    The solution of this problem involves two strong rarefaction waves and an intermediate state that is close to vacuum, what makes this problem  a hard test for numerical methods. ACAT methods  give stable solutions under CFL$\leq 1$ condition: Figure \ref{test61} shows the time evolution of the numerical results obtained with ACAT6. The smoothness indicators $\psi^3_{i+1/2}$ is also depicted: it can be seen how the discontinuities of  the first order derivatives are correctly captured. It can be also observed that, while at the rarefaction waves order 6 is selected,  lower accuracy is used at the constant regions close to the boundaries: this order reduction is due to the numerical oscillations produced by the 6th order method. A comparison of the different methods at time $t = 0.15$ is shown  in Figure \ref{test62} using 200-point mesh, where ACAT methods provide similar stable solutions. Although WENO solutions are stable, the third-order one is diffusive and the fifth-order one is oscillatory.

	\item{Right blast wave problem of Woodward $\&$ Colella:} the initial condition is

	\begin{equation} \label{t3incond}
	(\rho, u, p)
	= \left\{
	\begin{array}{ll}
	\displaystyle ( 1.0,0.0,1000) & \mbox {if } x < 1/2, \\
	\displaystyle (  1.0,0.0,0.01) & \mbox {if } x > 1/2.
	\end{array}\right.
	\end{equation}
	\begin{figure}[!ht]
		\setlength{\unitlength}{1mm}
		\centering	
		\includegraphics[width=\textwidth,height=6.9 cm]{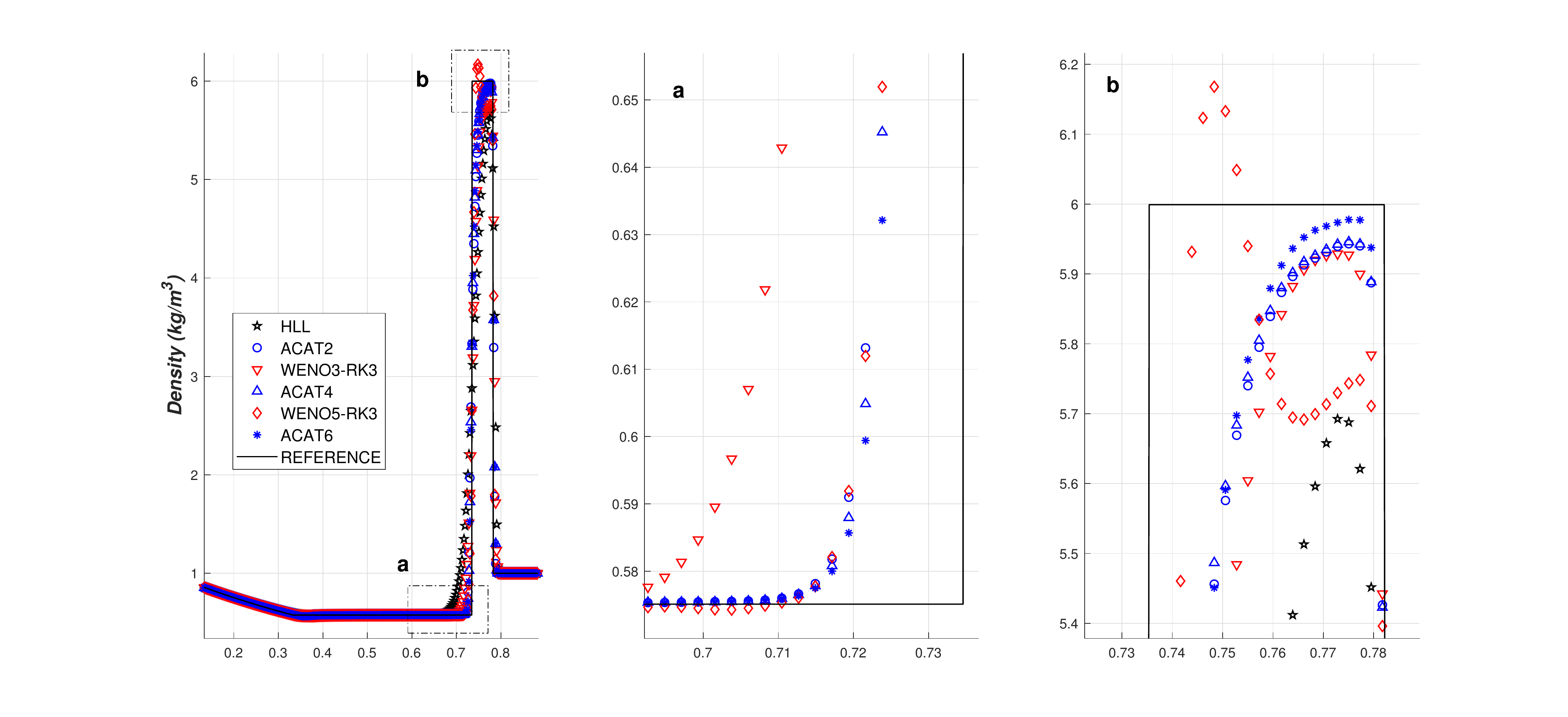}
		\vspace{-1.1 cm}
		\caption{ 1D Euler equations: right blast wave of the Woodward $\&$ Colella problem. Numerical solutions at time $t = 0.012$ using CFL$=0.8$ and $450$ points, (\textit{left}) and zooms close to  the shocks (\textit{center and right}).}
		\label{test7}
	\end{figure}
	
	 For this tests we use a 450-point mesh. The solution involves two strong shocks. Figure \ref{test7} shows the numerical densities obtained at time $t =0.012$: it can be observed that WENO methods produce oscillating solutions,  while ACAT methods give stable solutions whose accuracy increase with the order. In particular, this behavior can be seen in the two zooms close to the shocks.

Table \ref{table_1d_test_colella} shows the CPU time rates for this last one-dimensional test. 
A non-optimized implementation using  Matlab has been used for all the numerical methods. Therefore, this table has to be taken as a rough indication about computational cost. In particular, ACAT methods are highly parallelisable and do not need the storage of intermediate temporal stages: therefore, an optimized parallel implementation can lead to very different conclusions. With the implementations used here, ACAT2 is the cheapest method and its CPU time is taken as a reference. 
ACAT4 is competitive both in quality and computational cost compared to WENO-RK 3 and 5. The practical use of ACAT of order higher or equal than 6 requires an efficient implementation, otherwise the computational cost to increase the order is very big. The same happens with WENO-RK methods when the accuracy in time is increased due to the large number of stages required by SSPRK methods.

\begin{table}[htbp]
\begin{center}
\resizebox{8cm}{!}{
\begin{tabular}{|c|c|c|}
\hline
ACAT2     &     ACAT4     &     ACAT6  \\
1.00      &     5.88      &     12.46 \\
\hline
\hline
WENO3-RK3  &    WENO5-RK3     &       \\
2.86       &    5.08          &       \\
\hline 
\end{tabular}
}
\vspace{2mm}
\caption{CPU time rates for the Woodward and Colella problem.}
\label{table_1d_test_colella}
\end{center}
\end{table}

\end{itemize}

\subsection{2D Transport equation}
Let us consider the 2D transport equation
\begin{equation}\label{transport2d}
    u_t + a u_x + bu_y=0, \quad 
\end{equation}

with initial conditions 

\begin{equation}\label{ini_cond_trasnp2d}
u = \left\{
	\begin{array}{ll}
	1 & \mbox {if } x+y \leq 1/4, \\
	0 & \mbox {otherwise}.
	\end{array}\right.
\end{equation}

\begin{figure}[!ht]
	\setlength{\unitlength}{1mm}
	\centering	
	\begin{picture}(152,152)
	\put(4,3){\makebox(128,128)[c]{
			\includegraphics[width=\textwidth,height=10cm]{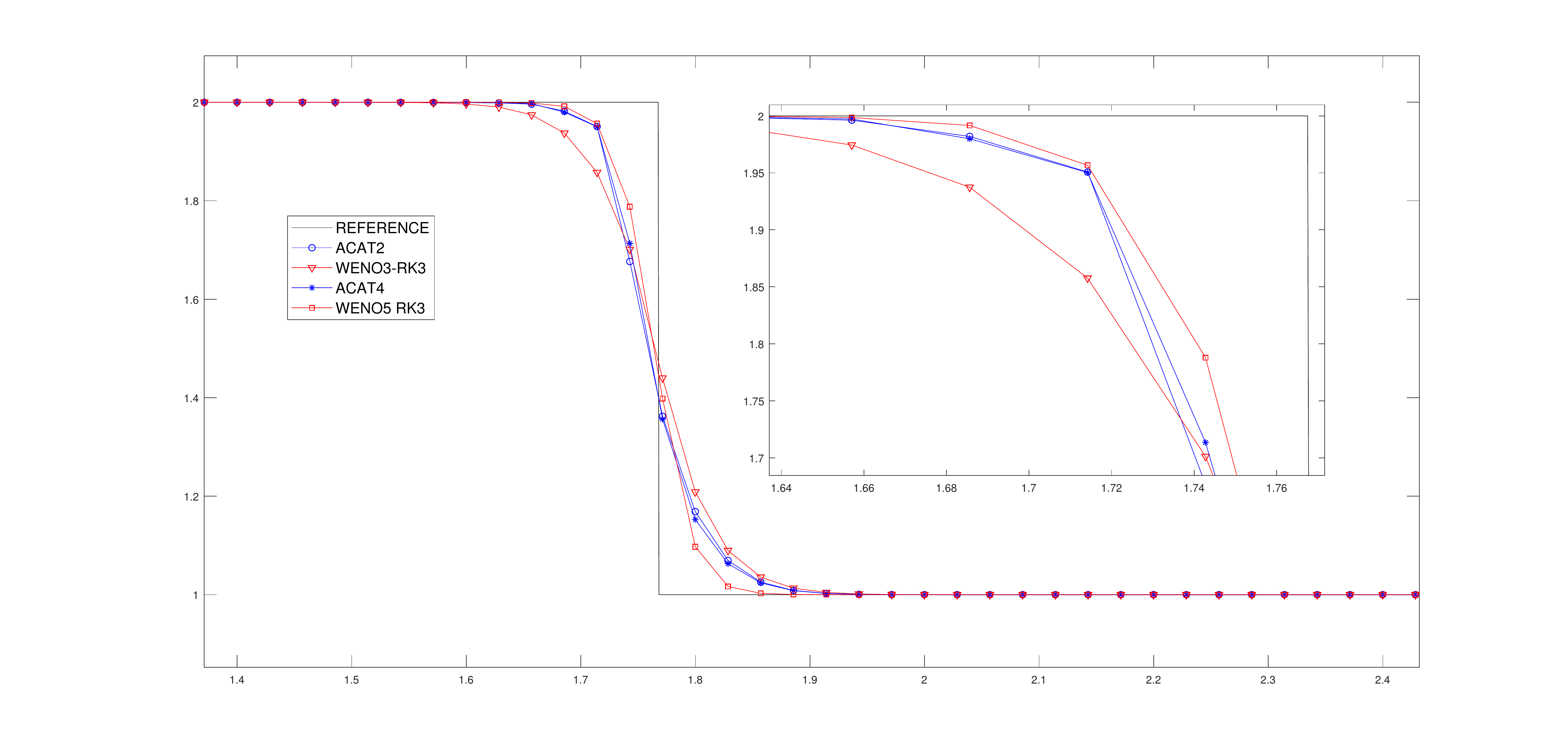}}}
	\end{picture}
	\vspace{-1.3 cm}
	\caption{2D Transport equation: test 1. Solution obtained with ACAT2, ACAT4, WENO3 RK3 and WENO5 RK3 at time $t = 1$: cut with a vertical plane passing through the line $y = x$ . Subplot: zoom close to the discontinuity}
	\label{transport_2d}
\end{figure}

We solve (\ref{transport2d}) on the spatial domain $[0,2]\times[0,2]$, using: $a,b=1$, $100\times 100$-point grid, CFL=0.5, free boundary conditions and $t=1$s. Figure \ref{transport_2d} shows a 1D cut over the line $y=x$ of the solutions obtained with  ACAT2, ACAT4, WENO3-RK3 and WENO5-RK3 at time $t = 1$.

\subsection{2D Euler equations}
Let us consider  the two-dimensional Euler equations for gas dynamics
\begin{equation}\label{2d_eq}
u_{t}+ {f}(u)_{x} + {g}(u)_{y}= 0, 
\end{equation}
where
$$
u=\left(\begin{array}{c}{\rho} \\ {\rho v} \\ {\rho w} \\ {E}\end{array}\right), 
\quad {f}(u)=\left(\begin{array}{c}{\rho v} \\ {\rho v^2 + p} \\ {\rho v w}  \\ {v(E+p)}\end{array}\right) , 
\quad {g}(u)=\left(\begin{array}{c}{\rho w} \\ {\rho v w} \\ {\rho w^2  + p} \\ {w(E+p)}\end{array}\right).
$$
Here, $\rho$ is the density; $v,w$ are the components of the velocity in the $x$ and $y$ directions; $E$,  the total energy per unit volume; $p$, the pressure. We consider the  equation of state 
\begin{equation}
p(\rho,v,w,E) = (\gamma - 1)(E-\frac{\rho}{2}(v^2+w^2)),
\end{equation}
and $\gamma$ is the ratio of specific heat capacities of the gas taken as 1.4.

We solve numerically (\ref{2d_eq}) using ACAT2 and ACAT4 for three of the nineteen configurations of the 2-D Riemann problems presented in \cite{Lax1998} whose initial conditions are given in Tables \ref{RP2d1}-\ref{RP2d2}.  These initial conditions consist of constant states at every quadrant of the spatial domain that are chosen so that the 1D Riemann problems corresponding to two adjacent states consist of only one one-dimensional simple wave: a shock $\textit{S}$, a rarefaction wave $\textit{R}$, or  a  slip line i.e. a contact discontinuity with discontinuous tangential velocity $\textit{J}$.  The sub-indexes  $(l,r) \in \{ (2,1),(3,2),(3,4),(4,1) \}$  indicate the involved quadrants. For shocks and rarefactions an over-arrow indicate the direction  (backward or forward). And for contact discontinuities a sign $+/-$ is used (instead of the over-arrow), to denote whether it is a positive or negative slip line.

These Riemann problems are numerically solved using a  $(400 \times 400)$-point grid and free boundary conditions. The CFL condition used to set the time steps is the following
\begin{equation*}
\Delta t = \frac{\mathrm{CFL}}{2} \,\, \min\left( \frac{\Delta x}{ \smax_x}, \frac{\Delta y}{\smax_y}\right),    
\end{equation*}
where 
$$ \smax_x = \max{i,j} \{ \abs{v_{i,j}^n} + c_{i,j}\}, \quad 
   \smax_y = \max{i,j} \{ \abs{w_{i,j}^n}+c_{i,j}\}, $$
with
$$ c = \sqrt{\frac{\gamma p}{\rho}}.$$
The CFL parameter is set to  $0.475$. 

\begin{table}[htbp]
\begin{center}
\resizebox{\textwidth}{!}{
\begin{tabular}{|llll|rrr|}
\hline
\textbf{Lax}  &  Configuration 4                &                &  \multicolumn{1}{c}{}  &   \multicolumn{3}{c|}{} \\ 
\hline
$p_2 = 0.35$    & $\rho_2 = 0.5065$    & $p_1 = 1.1$    & $\rho_1 = 1.1$                  &          &        &        \\
$u_2 = 0.8939$ & $v_2=0.0$      &$u_1= 0.0$    & $v_1=0.0$                  &          & $\overleftarrow{S}_{2,1}$      &        \\
$p_3=1.1$ & $\rho_3=1.1$ &  $p_4=0.35$   & $\rho_4=0.5065$               &  $\overrightarrow{S}_{3,2}$      &        & $\overrightarrow{S}_{4,1} $    \\ $u_3=0.8939$ & $v_3=0.8939$    &   $u_4=-0.0$   & $v_4=0.8939$                  &          & $\overleftarrow{S}_{3,4}$       &\\        
\hline
 \end{tabular}
}
\caption{2D Euler equations: test 1. Initial condition.}\label{RP2d1}
\end{center}
\end{table}

\begin{table}[htbp]
\begin{center}
\resizebox{\textwidth}{!}{
\begin{tabular}{|llll|rrr|}
\hline
\textbf{Lax }  &  Configuration 6                &                &  \multicolumn{1}{c}{}  &   \multicolumn{3}{c|}{} \\ 
\hline
$p_2 = 1.0$    & $\rho_2 = 2.0$     &$p_1 = 1.0$    & $\rho_1 = 1.0$                  &         &        &        \\
	 $u_2 = 0.75$ & $v_2=0.5$    & $u_1= 0.75$    & $v_1=-0.5$                  &         & $J^-_{2,1}$      &        \\
	 $p_3=1.0$ & $\rho_3=1.0$ &   $p_4=1.0$   & $\rho_4=3.0$               &  $J^+_{3,2}$      &        & $J^+_{4,1} $    \\
	 $u_3=-0.75$ & $v_3=0.5$      & $u_4=-0.75$   & $v_4=-0.5$                   &         & $J^-_{3,4}$       &  \\     
\hline
 \end{tabular}
}
\caption{2D Euler equations: test 2. Initial condition.}\label{RP2d2}
\end{center}
\end{table}

\begin{table}[htbp]
\begin{center}
\resizebox{\textwidth}{!}{
\begin{tabular}{|llll|rrr|}
\hline
\textbf{Lax }  &  Configuration 8                &                &  \multicolumn{1}{c}{}  &   \multicolumn{3}{c|}{} \\ 
\hline
$p_2 = 1.0$    & $\rho_2 = 1.0$    &  $p_1 = 0.4$    & $\rho_1 = 0.5197$                  &          &        &        \\
$u_2 = -0.6259$ & $v_2=0.1$     &$u_1= 0.1$    & $v_1=0.1$                  &         & $\overleftarrow{R}_{2,1}$      &        \\
$p_3=1.0$ & $\rho_3=0.8$   & $p_4=1.0$   & $\rho_4=1.0$               &  $J^-_{3,2}$      &        & $\overleftarrow{R}_{4,1} $    \\
$u_3=0.1$ & $v_3=0.1$     & $u_4=0.1$   & $v_4=-0.6259$                   &         & $J^-_{3,4}$       &\\      
\hline
 \end{tabular}
}
\caption{2D Euler equations: test 3. Initial condition.}\label{RP2d3}
\end{center}
\end{table}

\begin{figure}[!ht]
	\setlength{\unitlength}{1mm}
	\centering	
	\begin{picture}(152,180)
	\put(4,3){\makebox(128,170)[c]{
			\includegraphics[height=18cm]{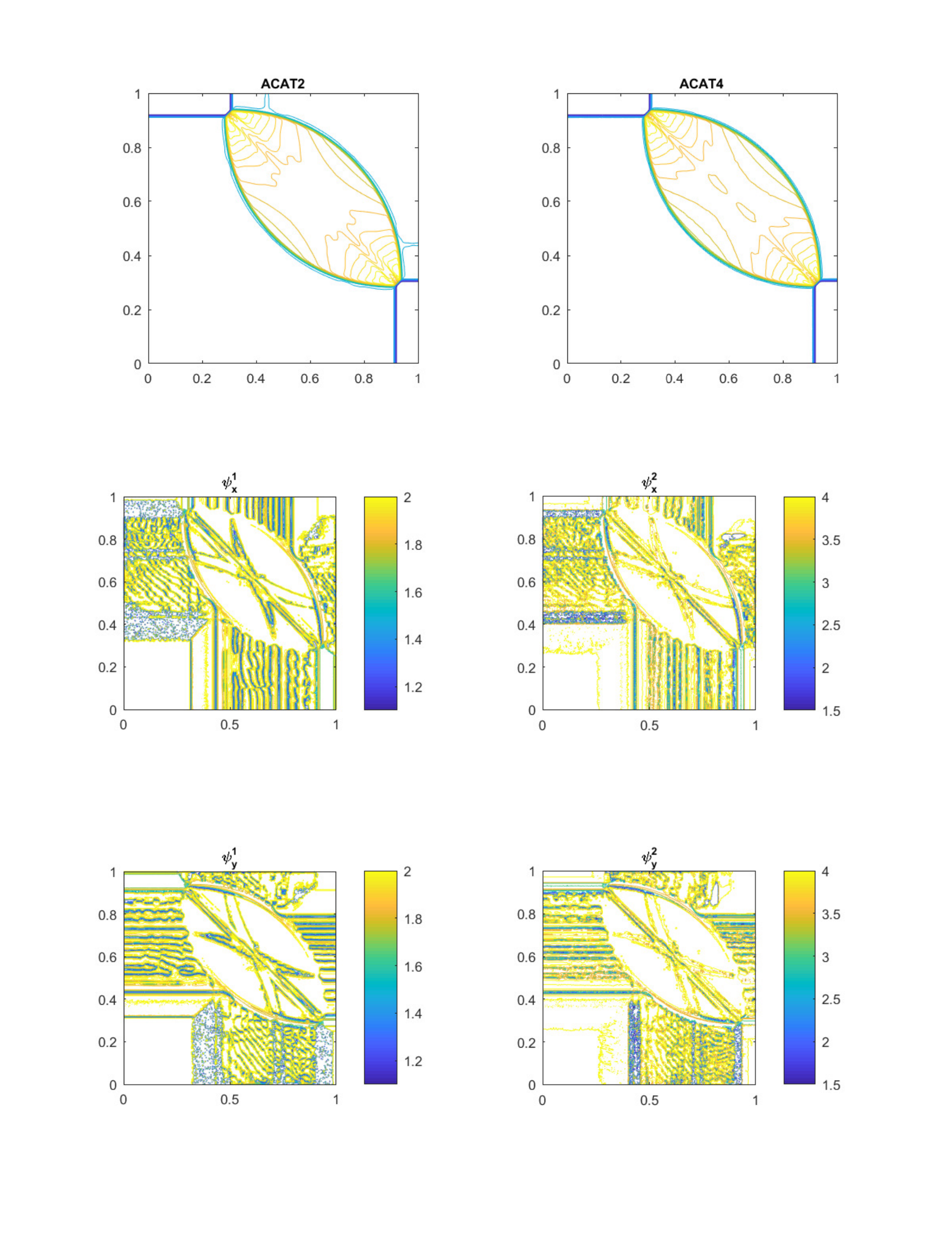}}}
	\end{picture}
	\vspace{-1 cm}
	\caption{2D Euler equations: test 1. Contour plots of the density at time $t = 0.25$ obtained with ACAT2 (\textit{left-up}) and ACAT4 ((\textit{right-up})). Contour plots of the smoothness indicators $\psi_x^1$(\textit{left-center}), $\psi_x^2$ (\textit{right-center}), $\psi_y^1$ (\textit{left-down}) and
	$\psi^2_y$ (\textit{right-down}).}
	\label{test8}
\end{figure}
\begin{figure}[!ht]
	\setlength{\unitlength}{1mm}
	\centering	
	\begin{picture}(152,180)
	\put(4,3){\makebox(128,170)[c]{
			\includegraphics[height=18 cm]{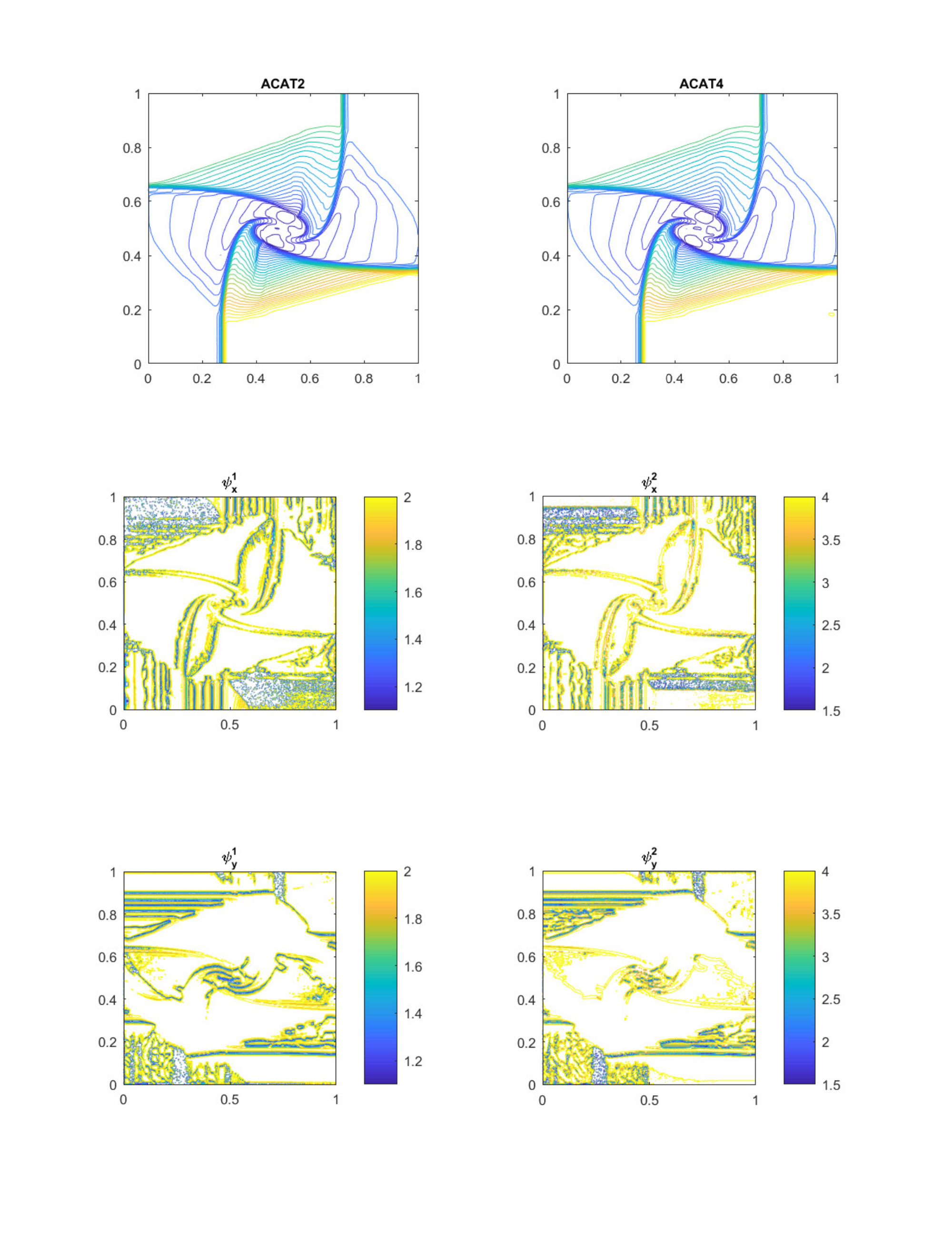}}}
	\end{picture}
	\vspace{-1 cm}
	\caption{2D Euler equations: test 2. Contour plots of the density at time $t = 0.3$ obtained with ACAT2 (\textit{left-up}) and ACAT4 (\textit{right-up}). Contour plots of the smoothness indicators $\psi_x^1$ (\textit{left-center}), $\psi_x^2$ (\textit{right-center}), $\psi_y^1$ (\textit{left-down}) and
	$\psi^2_y$ (\textit{right-down}).}
	\label{test9}
\end{figure}
\begin{figure}[!ht]
	\setlength{\unitlength}{1mm}
	\centering	
	\begin{picture}(152,180)
	\put(4,3){\makebox(128,170)[c]{
			\includegraphics[height=18cm]{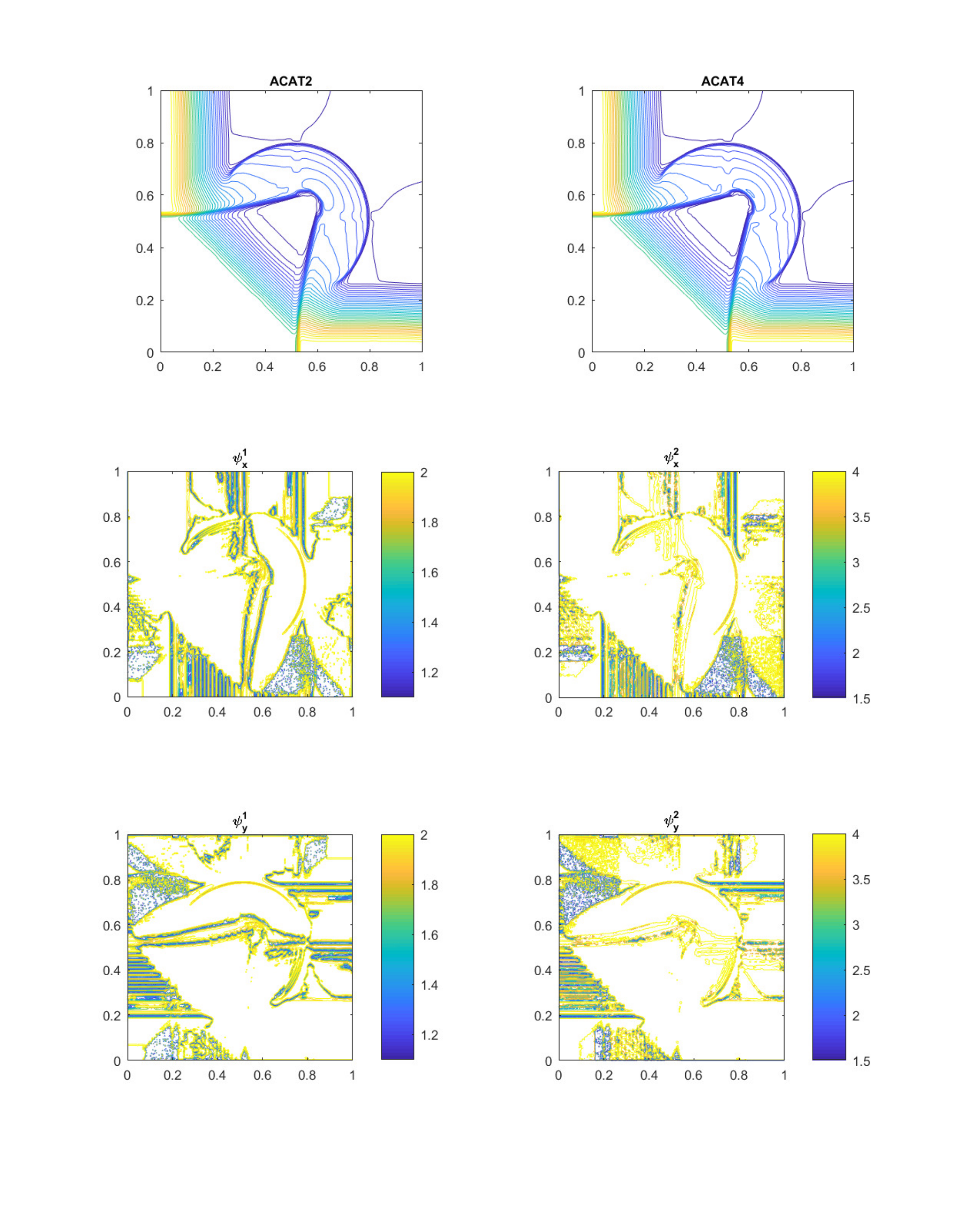}}}
	\end{picture}
	\vspace{- 1 cm}
	\caption{2D Euler equations: test 3. Contour plots of the density at time $t = 0.25$ obtained with ACAT2 (\textit{left-up}) and ACAT4 (\textit{right-up}). Contour plots of the smoothness indicators $\psi_x^1$(\textit{left-center}), $\psi_x^2$ (right-center), $\psi_y^1$ (\textit{left-down}) and
	$\psi^2_y$ (\textit{right-down}).}
	\label{test10}
\end{figure}

\begin{figure}[!ht]
	\setlength{\unitlength}{1mm}
	\centering	
	\begin{picture}(152,152)
	\put(4,3){\makebox(128,128)[c]{
			\includegraphics[width = 12.4cm,height=12cm]{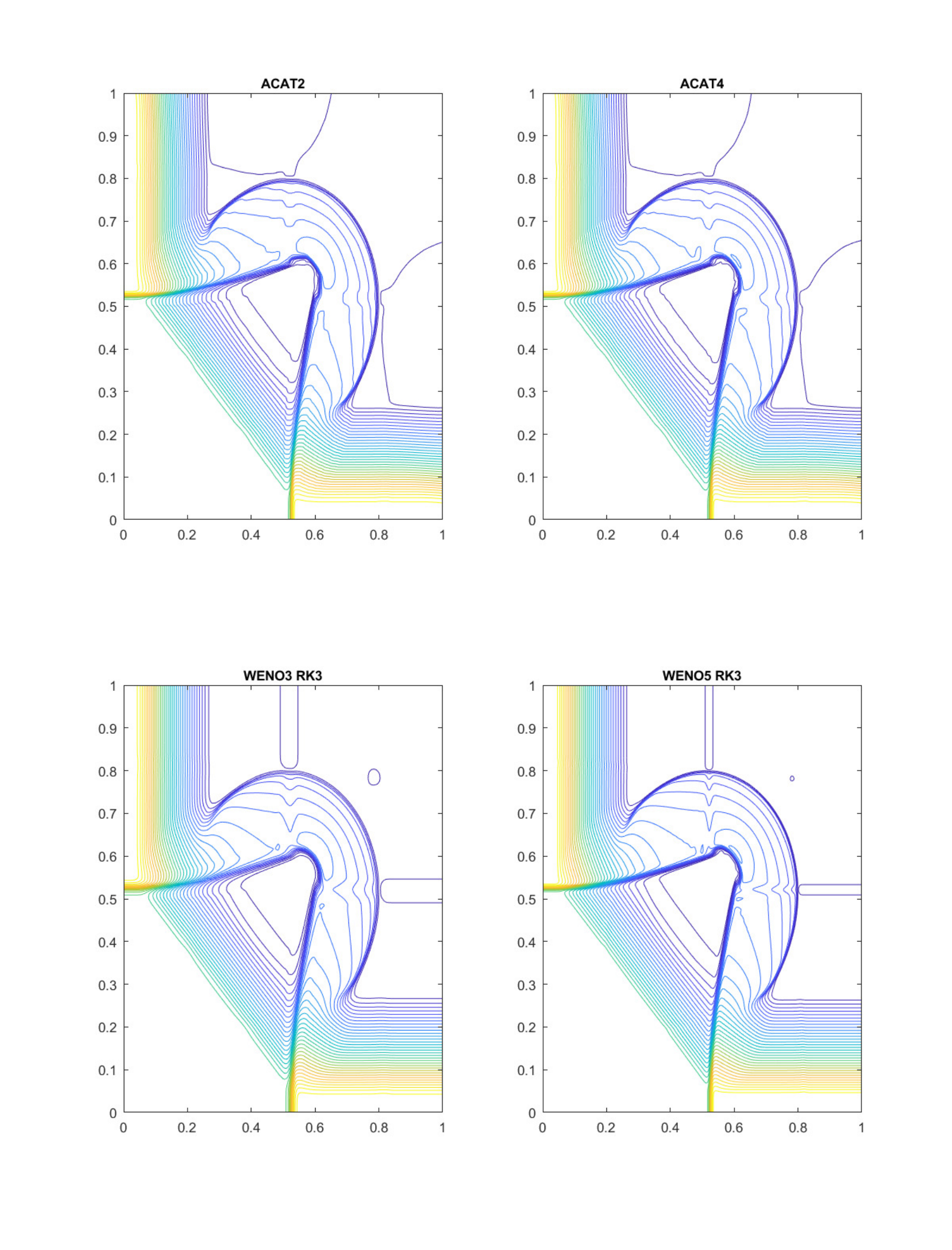}}}
	\end{picture}
	\vspace{-1 cm}
	\caption{2D Euler equations: test 3. Contour plots of the density at time $t = 0.25$ obtained with ACAT2 (\textit{left-up}), ACAT4 (\textit{right-up}), WENO3 RK3 (\textit{left-down}) and WENO5 RK3 (\textit{right-down}).}
	\label{test10_2}
\end{figure}

Figures \ref{test8}, \ref{test9} and \ref{test10} show the numerical solutions for the density given by ACAT2 and ACAT4. We include in each figure a general view of the numerical density given by ACAT2  (left-up) and ACAT4 (right-up);  the smoothness indicators $\psi^1_x$ (left-center) and $\psi^2_x$ (right-center) in  the $x$-direction; the smoothness indicators $\psi^1_y$ (left-down) and $\psi^2_y$ (right-down) in  the $y$-direction. In all cases, the solutions are stable and  similar to those obtained in \cite{KT2002} with a finite volume method. Observe how the indicators $\psi_x^2$ and $\psi_y^2$ detect better the smoothness regions than $\psi_x^1$ and $\psi_y^1$, what implies a better resolution in the numerical solutions obtained with ACAT4. However, the computational cost increases with the order as it happened for 1d problems, see Table \ref{table_2d_tiempos}.

\begin{table}[htbp]
\begin{center}
\resizebox{8cm}{!}{
\begin{tabular}{|c|c|c|c|c|}
\hline
ACAT2     &     ACAT4     &     ACAT6 &  WENO3-RK3  &    WENO5-RK3    \\
\hline
1.00      &     9.98      &     96.91  & 3.23       &    9.968  \\
\hline 
\end{tabular}
}
\vspace{2mm}
\caption{2D Euler equations test 3: CPU time rates. }
\label{table_2d_tiempos}
\end{center}
\end{table}

\noindent In Figure \ref{test10_2} the numerical densities obtained with  ACAT2, ACAT4, WENO3 RK3, and WENO5 RK5 at time $ t = 0.25$ are compared. 

\newpage
\section{Conclusions} \label{sec4}

In this work, the new family of high-order shock-capturing Adaptive Compact Approximate Taylor (ACAT) methods for systems of conservation laws has been introduced. These method are an order adaptive version of the Compact Approximate Taylor Methods introduced in \cite{CP2019} in which the solution at every point is updated using the stencil of maximal length for which the solution is smooth.

The 5-point stencil ACAT2 method coincides with the FL-CAT2 introduced in \cite{CP2019} that combines CAT2 with a first order robust numerical method using a standard flux limiter. For higher orders, a new family of smoothness indicators has been introduced to select the maximal length stencil at every point.
The expression of the methods for 1D or 2D systems of conservation laws has been given. 

The results obtained with the new family of methods in a number of test cases have been compared  with the corresponding WENO-RK method (Finite 
Differences WENO reconstructions in space, TVD-RK in time). The linear transport equation, Burgers equation, the 1D
and 2D compressible Euler equations have been considered. For $\mathrm{CFL} \le 0.5$ all the numerical methods work correctly, 
and the results obtained with WENO or ACAT methods are similar. Nevertheless, for
CFL values  close to one, ACAT still give good results while WENO methods may be oscillatory. The possibility of using larger time steps compensate the extra 
computational cost of a temporal iteration. 
ACAT methods are more 
expensive in computational time and number of operations due to its local character.
Nevertheless, with the non-optimized implementation of the methods performed to solve the test cases shown here, the computational cost to increase the order from 4 to 6  is very big, specially for 2D problems: an optimized implementation is necessary to exploit all the potentialities of these methods that are highly parallelisable and do not need the storage of intermediate temporal stages.
Further developments include:
\begin{itemize}
    \item An optimized implementation in GPU architectures.
    \item The extension to systems of balance laws.
    
\end{itemize}

\section*{Acknowledgements} 
 This research has received funding from the European Union’s Horizon 2020 research and innovation program, under the Marie Sklodowska-Curie grant agreement No 642768. It  has been also partially supported by the Spanish Government and FEDER through the Research project RTI2018-096064-B-C21. E. Macca has been also partially supported by the Piano triennale della Ricerca 2016-2018, Department of Mathematics and Computer Sciences, University of Catania. D. Zor\'io is also supported by Fondecyt Project 3170077.

\section*{Appendix}

\begin{figure}[!ht]
	\setlength{\unitlength}{1mm}
	\centering	
	\includegraphics[width=\textwidth,height=16
	cm]{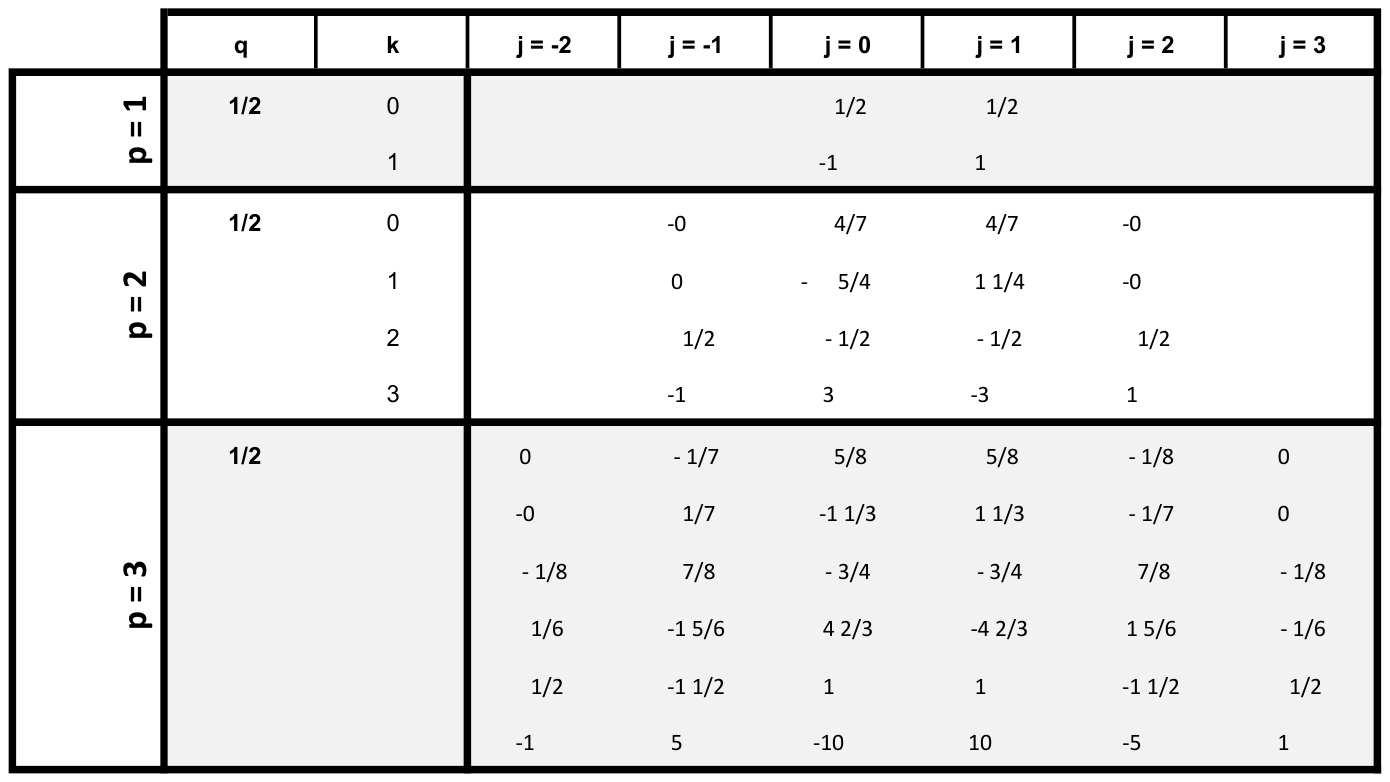}
	\vspace{-7 cm}
	\caption{ The $\delta^k_{p,j}$ coefficients of the differentiation formula \eqref{F}  for $p=1,2,3$. }
	\label{anexo1}
\end{figure}

\begin{figure}[!ht]
	\setlength{\unitlength}{1mm}
	\centering	
	\includegraphics[width=\textwidth,height=16
	cm]{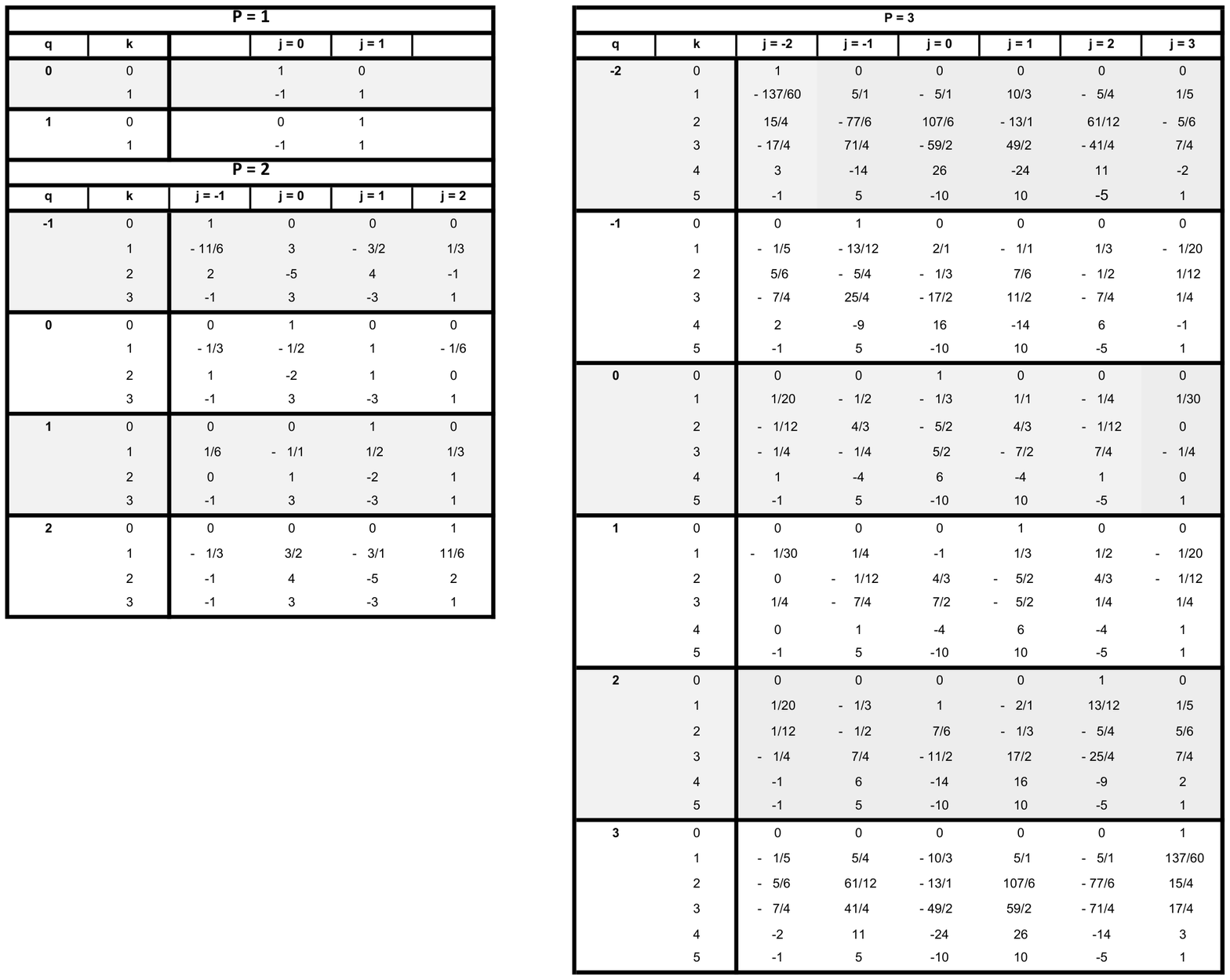}
	\vspace{-2 cm}
	\caption{The $\gamma^{k,q}_{p,j}$ coefficients of the differentiation formula \eqref{upwF}  for $p=1,2,3$.}
	\label{anexo2}
\end{figure}

The coefficients $\delta^k_{p,j}$ and $\gamma^{k,q}_{p,j}$ of the differentiation formulas \eqref{F} and \eqref{upwF} for $p=1,2,3$ are shown in Figures \ref{anexo1} and \ref{anexo2} respectively. Algorithms to compute those coefficients can be found in \cite{Fornberg1} and \cite{CP2019}.

\newpage

\bibliographystyle{unsrt}
\bibliography{mybibfile}

\begin{thebibliography}{10}

\bibitem{CP2019}
H.~Carrillo and C.~Parés.
\newblock Compact approximate {T}aylor methods for systems of conservation
  laws.
\newblock {\em Journal of Scientific Computing}, 80:1832--1866, 2019.

\bibitem{ZBM2017}
D.~Zor{\'i}o, A.~Baeza, and P.~Mulet.
\newblock An approximate lax--wendroff-type procedure for high order accurate
  schemes for hyperbolic conservation laws.
\newblock {\em Journal of Scientific Computing}, 71:246--273, 2017.

\bibitem{LeVeque2007book}
R.J. LeVeque.
\newblock {\em Finite difference methods for ordinary and partial differential
  equations: steady-state and time-dependent problems (Classics in Applied
  Mathematics)}.
\newblock Society for Industrial and Applied Mathematics, Philadelpia, PA.
  USA., 1 edition, 2007.

\bibitem{Toro2009book}
E.F. Toro.
\newblock {\em {Riemann} Solvers and Numerical Methods for Fluid Dynamics}.
\newblock Springer, third edition, 2009.

\bibitem{Gideon1971}
G.~Zwas and S.~Abarbanel.
\newblock Third and fourth order accurate schemes for hyperbolic equations of
  conservation law form.
\newblock {\em Mathematics of Computation}, 25(114):229--236, 1971.

\bibitem{Ader2001}
E.F. Toro, R.C. Millington, and L.A.M Nejad.
\newblock Towards very high order godunov schemes.
\newblock {\em Godunov Methods. Theory and Applications E.F. Toro ed.,
  Kluwer/Plenum Academic Publishers}, pages 907--940, 2001.

\bibitem{TitarevToro2002}
V.A. Titarev and E.F. Toro.
\newblock {ADER}: Arbitrary high order {G}odunov approach.
\newblock {\em Journal of Scientific Computing}, 17:609--618, 2002.

\bibitem{Schwartzkopff2002}
T.~Schwartzkopff, C.~D. Munz, and E.F. Toro.
\newblock A high-order approach for linear hyperbolic systems in 2d.
\newblock {\em Journal of Scientific Computing}, 17:231--240, 2002.

\bibitem{DumbserToro2008}
C.~Enaux, M.~Dumbser, and E.F. Toro.
\newblock Finite volume schemes of very high order of accuracy for stiff
  hyperbolic balance laws.
\newblock {\em Journal of Computational Physics}, 227(2):3971--4001, 2008.

\bibitem{PNPM}
M.~Dumbser, D.~Balsara, E.F. Toro, and C.D. Munz.
\newblock A unified framework for the construction of one-step finite-volume
  and discontinuous galerkin schemes.
\newblock {\em Journal of Computational Physics}, 227:8209--8253, 2008.

\bibitem{CPZ2020}
H.~Carrillo, C.~Parés, and D.~Zorío.
\newblock Approximate {T}aylor methods with fast and optimized weighted
  essentially non-oscillatory reconstructions.
\newblock {\em arXiv:2002.08426v1 [math.NA] 19 Feb 2020}, 2020.

\bibitem{Fornberg1}
B.~Fornberg.
\newblock Generation of finite difference formulas on arbitrarily spaced grids.
\newblock {\em Mathematics of Computation}, 51:699--706, 1988.

\bibitem{ZBBM2019}
A.~Baeza, R.~Bürger, P.~Mulet, and D.~Zorío.
\newblock On the efficient computation of smoothness indicators for a class of
  weno reconstructions.
\newblock {\em Journal of Scientific Computing}, 80:1240–1263, 2019.

\bibitem{BBMZO3}
A.~Baeza, R.~Bürger, P.~Mulet, and D.~Zorío.
\newblock An efficient third-order {WENO} scheme with unconditionally optimal
  accuracy.
\newblock {\em SIAM Journal on Scientific Computing (To appear)}, 2020.

\bibitem{Qiu2002}
J.~Qiu and C.-W. Shu.
\newblock On the construction, comparison, and local characteristic
  decomposition for high-order central weno schemes.
\newblock {\em J. Comput. Phys.}, 183(1):187–209, 2002.

\bibitem{Roesuperbee}
P.L. Roe.
\newblock Characteristic-based schemes for the euler equations.
\newblock {\em Annu. Rev. Fluid Mech.}, 18:337--365, 1986.

\bibitem{Shu1997}
C.~W. Shu.
\newblock Essentially non-oscillatory and weighted essentially non--oscillatory
  schemes for hyperbolic conservation laws.
\newblock Technical report, Institute for Computer Applications in Science and
  Engineering (ICASE), 1997.

\bibitem{Gottlieb2011}
S.~Gottlieb, D.~Ketcheson, and C.W. Shu.
\newblock {\em Strong Stability Preserving Runge-Kutta and multistep time
  discretizations}.
\newblock Word Scientific, 1 edition, 2011.

\bibitem{Sod1978}
G.A. Sod.
\newblock A survey of several finite difference methods for systems of
  nonlinear hyperbolic conservation laws.
\newblock {\em Journal of Computational Physics}, 27(1):1--31, 1978.

\bibitem{Einfeldt1991}
B.~Einfeldt, P.L Roe, C.D. Munz, and B.~Sjogreen.
\newblock {On Godunov--type methods near low densities}.
\newblock {\em Journal of Computational Physics}, 92:273--295, feb 1991.

\bibitem{Wood1984}
P.~Woodward and P.~Colella.
\newblock The numerical simulation of two-dimensional fluid flow with strong
  shocks.
\newblock {\em Journal of Computational Physics}, 1:115--173, 1984.

\bibitem{Lax1998}
P.~Lax and Liu Xu-Dong.
\newblock Solution of two-dimensional riemann problems of gas dynamics by
  positive schemes.
\newblock {\em SIAM Journal on Scientific Computing}, 19F(2):319--340, 1998.

\bibitem{KT2002}
A.~Kurganov and E.~Tadmor.
\newblock Solution of two-dimensional riemann problems for a gas dynamics
  without riemann problem solvers.
\newblock {\em Numer. Methods Partial Differential Equations}, 18:584--608,
  2002.

\end{thebibliography}

\end{document}